\documentclass[aap,MSNbibl,seceqn,dvips]{arximspdf}
\usepackage{graphicx}
\doi{10.1214/13-AAP927} 
\volume{24}
\issue{1}
\pubyear{2014}
\firstpage{378}
\lastpage{421}

\makeatletter

\newcommand{\GI}{\mathit{GI}}

\newcommand{\rrvert}{\vert}
\newcommand{\llvert}{\vert}

\newcommand{\cal}{\mathcal}

\newtheorem{theorem}{Theorem}[section]
\newtheorem{corollary}{Corollary}[section]

\newproclaim{remark}{Remark}[section]

\newcommand{\RR}{{\mathbb R}}
\newcommand{\sB}{{\cal B}}
\newcommand{\sU}{{\cal U}}
\newcommand{\sR}{{\cal R}}
\newcommand{\sK}{{\cal K}}
\newcommand{\ra}{\rightarrow}
\newcommand{\deq}{\stackrel{\mathrm{d}}{=}}
\newcommand{\D}{\mathbb{D}}
\newcommand{\C}{\mathbb{C}}
\newcommand{\ep}{\varepsilon}

\makeatother

\begin{document}
\begin{frontmatter}

\title{Many-server heavy-traffic limit for queues
with time-varying parameters\thanksref{T1}}
\pdftitle{Many-server heavy-traffic limit for queues
with time-varying parameters}
\runtitle{Many-server heavy-traffic limits}

\thankstext{T1}{Supported by NSF Grants CMMI 0948190 and 1066372.}

\begin{aug}
\author[A]{\fnms{Yunan} \snm{Liu}\ead[label=e1]{yliu48@ncsu.edu}}
\and
\author[B]{\fnms{Ward} \snm{Whitt}\corref{}\ead[label=e2]{ww2040@columbia.edu}}
\runauthor{Y. Liu and W. Whitt}
\affiliation{North Carolina State University and Columbia University}
\address[A]{Department of Industrial\\
\quad and Systems Engineering \\
North Carolina State University \\
Room 446, 400 Daniels Hall \\
Raleigh, North Carolina 27695 \\
USA\\
\printead{e1}}
\address[B]{Department of Industrial Engineering\\
\quad and Operations Research\\
Columbia University\\
New York, New York 10027-6699 \\
USA\\
\printead{e2}}
\end{aug}

\received{\smonth{9} \syear{2011}}
\revised{\smonth{2} \syear{2013}}

%
\begin{abstract}
A many-server heavy-traffic FCLT is proved for the $G_t/M/s_t+\GI$
queueing model, having time-varying arrival rate and staffing, a
general arrival process satisfying a FCLT, exponential service times
and customer abandonment according to a general probability
distribution. The FCLT provides theoretical support for the
approximating deterministic fluid model the authors analyzed in a
previous paper and a refined Gaussian process approximation, using
variance formulas given here. The model is assumed to alternate between
underloaded and overloaded intervals, with critical loading only at the
isolated switching points. The proof is based on a recursive analysis
of the system over these successive intervals, drawing heavily on
previous results for infinite-server models. The FCLT requires careful
treatment of the initial conditions for each interval.
\end{abstract}

%
\begin{keyword}[class=AMS]
\kwd[Primary ]{60K25}
\kwd[; secondary ]{60F17}
\kwd{90B22}
\end{keyword}
\begin{keyword}
\kwd{Many-server queues}
\kwd{queues with time-varying arrivals}
\kwd{nonstationary queues}
\kwd{customer abandonment}
\kwd{nonexponential patience distribution}
\kwd{heavy traffic}
\kwd{functional central limit theorem}
\kwd{Gaussian approximation}
\kwd{deterministic fluid approximation}
\end{keyword}

\pdfkeywords{60K25, 60F17, 90B22, Many-server queues,
queues with time-varying arrivals,
nonstationary queues, customer abandonment,
nonexponential patience distribution,
heavy traffic, functional central limit theorem,
Gaussian approximation,
deterministic fluid approximation}

\end{frontmatter}

\section{Introduction}\label{secIntro}

This paper is a sequel to \cite{LW10}, in which we developed and
analyzed a deterministic fluid model approximating the
$G_t/\GI/s_t+\GI$ queueing model, having a general arrival process with
time-varying arrival rate (the initial~$G_t$),
independent and identically distributed (i.i.d.)
service times with a general cumulative distribution function (c.d.f.) $G$
(the~first~$\GI$),
a time-varying large number of servers (the~$s_t$)
and customer abandonment from queue
with i.i.d. patience times with a general c.d.f. $F$ (the~final $+\GI$).
The fluid model was assumed to alternate between intervals of
underloading (UL) and overloading (OL).
We conducted simulation experiments showing that the fluid
approximation is effective for approximating individual sample paths of
stochastic processes
of very large systems (e.g., with hundreds of servers) and the mean
values of smaller systems (e.g., with tens of servers,
provided that these systems are not critically loaded or too nearly so).
See \cite{FMMW08,GKW07,HM10,MMR98,M02,N82} for background on methods
to analyze the performance of queues with time-varying arrival rates
and their application.

The present paper establishes many-server heavy-traffic limits
that provide mathematical support for both the previous fluid
approximation and a refined Gaussian process approximation in the special
case of exponential ($M$) service times.
Based directly on the limit theorems here, we propose approximating the
time-varying number of the customers in system (including those in
service and those waiting in queue if any), $X_n (t)$,
by a Gaussian distribution for each $t$, in particular,
%
%
\begin{equation}
\label{Xapprox} X_n (t) \approx n X (t) + \sqrt{n} \hat{X}(t) \deq N
\bigl(n X (t), n\sigma ^2_{\hat{X}} (t)\bigr),
\end{equation}
where $N(m,\sigma^2)$ denotes a Gaussian random variable with mean $m$
and variance~$\sigma^2$,
$X (t)$ is the deterministic fluid approximation proposed and analyzed
previously in \cite{LW10},
and now supported by the functional weak laws of large numbers
(FWLLNs) in Theorems \ref{thFWLLN} and \ref{thUL}, while
$\hat{X} (t)$ is a zero-mean Gaussian process with variance $\sigma
^2_{\hat{X}} (t) \equiv \operatorname{Var}(\hat{X} (t))$
obtained\vspace*{2pt} from the functional central limit theorems (FCLTs) in
Theorems \ref{thFCLT}, \ref{thOther} and \ref{thUL}.
Explicit formulas for the variance function $\sigma^2_{\hat{X}} (t)$
are given in Corollary \ref{corVar} to go\vspace*{1pt} with the
explicit expressions for the fluid function $X (t)$ determined
previously in \cite{LW10}, and reviewed here in Section~\ref{secFluid}.

As in \cite{LW10}, we assume that the system alternates between UL
intervals and OL intervals, where the system loading is determined by
the fluid model, which has the same parameters; that is, the system is
said to be UL (OL) if the
fluid model is UL (OL).
Sufficient conditions for the fluid model to alternate between OL and
UL intervals were given in Section~3 of \cite{LW11a}.
In the terminology of many-server heavy-traffic limits \cite{GMR02},
that means that the system alternates between
quality-driven (QD) UL regimes and efficiency-driven (ED) OL regimes.
We assume that the
system is never critically loaded, that is, in the
quality-and-efficiency-driven (QED) regime, except at the isolated regime
switching points. That allows us to apply previous results for infinite-server
queues in \cite{PW10} in our analysis of both UL and OL intervals.

Explicitly avoiding the QED regime runs counter to
most of the extensive research on many-server queues, for example, as
in \cite{GMR02,HW81,KR10,KR11,R09}.
However, we think the alternating UL and OL model can provide useful
approximations because
it provides mathematical simplification. This regime has engineering
relevance with time-varying arrivals
because many systems are unable to dynamically adjust staffing
to respond adequately to time-varying arrival rates, and thus do
experience periods of overloading.
Hospital emergency rooms are examples.\vadjust{\goodbreak}

The limits here extend
previous limits for the Markovian $M_t/M/s_t+M$ models with
time-varying arrival rates and staffing in \cite{MMR98,MMRR99,MMRS99,P08}.
To treat the model with general patience distribution, we exploit
limits for two-parameter stochastic processes in
infinite-server models in \cite{PW10}; also see \cite{KP97,RT09}.
Heavy-traffic limits for the stationary $G/M/s+\GI$ model were
established in \cite{ZM05},
where references on previous work can be found.
A previous discrete-time many-server limit for the $G_t/\GI/s+\GI$ model
with time-varying arrivals was established in \cite{W06}; in contrast,
here the limit is for a model with smooth parameters.
In a sequel to this paper, \cite{LW12}, we establish a FWLLN for the
more general $G_{t}/\GI/s_{t}+\GI$ model.
It remains to extend the FCLT to nonexponential service times.

In \cite{LW10} we saw that the analysis of the performance of the
$G_t/M/s_t +\GI$ fluid model depends critically on a careful analysis of the
waiting time
of the fluid at the head of the line (that has been waiting in queue
the longest).
That fluid head-of-the-line waiting time (HWT) $w(t)$ was identified by
carefully relating the new service capacity becoming available
due to service completion and changing capacity to the flow into
service from the queue.
That led to an ordinary differential equation (ODE) characterizing the
deterministic HWT function $w(t)$,
proved in Theorem~3 of \cite{LW10} and reviewed here in (\ref{wODE}).
Closely paralleling that ODE, we find that the stochastic limit process
for the FCLT-scaled HWT, $\hat{W} (t)$,
is characterized by a stochastic differential equation (SDE); see (\ref{sde}).

We primarily focus on the number in system $X_n (t)$, as in (\ref{Xapprox}),
because that process and the associated FCLT-scaled version [see (\ref
{fcltScaled}) below]
tends to be better behaved than the number in queue, $Q_n (t)$, and the
number in service,
$B_n (t)$, and the associated FCLT-scaled versions of them.
This is reflected by the limit processes for the FCLT-scaled versions.
For each $t$ in the interior of an OL interval,
$(\hat{Q} (t), \hat{B} (t)) = (\hat{X} (t), 0)$;
for each $t$ in the interior of an UL interval,
$(\hat{Q} (t), \hat{B} (t)) = (0, \hat{X} (t))$;
for each switching point $t$, $(\hat{Q} (t), \hat{B} (t)) = (\hat{X}
(t)^+, \hat{X} (t)^-)$,
where $(x)^+ \equiv\max{\{x, 0\}}$ and $(x)^- \equiv\min{\{x, 0\}}$.
Thus,\vspace*{1pt} in contrast to $\hat{X}$, which has continuous sample paths, the
sample paths of $\hat{Q}$ and $\hat{B}$
are discontinuous and are typically neither right-continuous nor
left-continuous at each switching point.

Thus, even though limits can be obtained
for FCLT-scaled versions of the number in queue, $Q_n (t)$, and the
number in service, $B_n (t)$,
yielding approximations such as $Q_n (t) \approx n Q (t) + \sqrt{n}
\hat{Q}(t)$, paralleling (\ref{Xapprox}),
we instead suggest approximating these processes by truncating the
number in system $X_n (t)$ with respect to the time-varying service
capacity $s_n (t)$;
that is, we propose the alternative approximations
%
%
\begin{eqnarray}
\label{QBapprox} Q_n (t) & = & \bigl(X_n (t) -
s_n (t)\bigr)^{+} \approx\bigl(n \bar{X} (t) + \sqrt {n}
\hat{X}(t) - s_n (t)\bigr)^{+},
\nonumber\\[-8pt]\\[-8pt]
B_n (t) & = & X_n (t) \wedge s_n (t)
\approx\bigl(n \bar{X} (t) + \sqrt {n} \hat{X}(t)\bigr) \wedge s_n
(t),\nonumber
\end{eqnarray}
exploiting (\ref{Xapprox}). This approximation is
convenient because formulas are known for the means and variances of
such truncated Gaussian variables.

We study such refined engineering approximations based on the
many-server heavy-traffic limits established here, including (\ref{QBapprox}),
in a future paper \cite{LW11b}.
However, immediate insight can be obtained by considering the special
case of the $M_t/M/s_t+M$ model with abandonment rate $\theta$ equal
to the service rate $\mu$.
As discussed in Section~6 of \cite{FMMW08}, the number in system in
this model is distributed the same as in the associated
$M_t/M/\infty$ model, for which the number in system at each time has
a Poisson distribution. In this case,
the approximation by (\ref{Xapprox}) is known to perform very well.
In that context, clearly the approximations in (\ref{QBapprox}) perform
well too, whereas even the direct approximation for the means, $n Q(t)$
and $n B(t)$ do not perform well near critical loading.

Here is how the rest of this paper is organized:
in Section~\ref{secSeq} we specify the sequence of $G_t/M/s_t+\GI$
queueing models we consider
and the associated scaled stochastic processes for the FWLLN and the FCLT.
In Section~\ref{secFluid} we review the $G_t/M/s_t+\GI$ fluid model,
which arises as the limit in the FWLLN
and provides centering terms for the FCLT.
In Section~\ref{secOL} we state the new results for each OL interval,
while in Section~\ref{secUL} we state the (easier) new results
for each UL interval. In Section~\ref{secOLproofs} we prove the FWLLN
and FCLT for OL intervals;
in Section~\ref{secCors} we prove two corollaries for OL intervals; finally,
in Section~\ref{secULproofs} we prove the FWLLN and FCLT for UL
intervals. In order to confirm that the formulas for the
variances given in Corollary \ref{corVar} are correct, thus providing
practical confirmation for all the results, we conduct simulation
experiments of both large and small queueing systems in Section~\ref{secCompare}.
We conclude in Section~\ref{secRefined} by discussing an extension
with extra $\sqrt{n}$ terms in the arrival rates and staffing functions.

\section{A sequence of $G_t/M/s_t+\mathit{GI}$ models}\label{secSeq}

In this paper we consider a sequence of $G_t/M/s_t+\GI$ queueing models
indexed by $n$.
Model $n$ has a general arrival process with time-varying arrival rate
$\lambda_n (t) \equiv n \lambda(t)$, i.i.d. exponential service times
with cumulative distribution function (c.d.f.) $G (t) \equiv1 -e^{- \mu
t}$, a time-varying number of servers $s_n (t) \equiv\lceil n
s(t)\rceil$
[the least integer above $n s(t)$]
and customer abandonment from queue, where
the patience times of successive customers to enter queue are i.i.d.
with general c.d.f. $F$, where we assume that $F$ is differentiable,
with probability density function (p.d.f.) $f$
with $F^c (x) > 0$ and $f (x) > 0$ for all $x$. Our scaling of the
fixed functions $\lambda$ and $s$
induces the familiar many-server heavy-traffic scaling; the functions
$\lambda$ and $s$ are the arrival rate and staffing level in the
associated fluid model,
assumed to be suitably smooth,
as specified in the next section.
The arrival process, service times and patience times are mutually independent.
New arrivals enter service immediately if there is a free server;
otherwise they join the queue, from which they enter service in order
of arrival,
if they do not first abandon.

Let $\D\equiv D(I)$ be the usual space of right-continuous real-valued
functions with left limits on a subinterval $I$ of $\RR$,
endowed with the Skorohod $J_1$ topology, which for continuous limits
reduces to uniform convergence over all compact subintervals of $I$.
Let $\Rightarrow$ denote convergence in distribution \cite{W02}.
Let $N_n (t)$ count the number of arrivals in $[0,t]$.
We assume that the sequence of arrival processes $\{N_n\}$ satisfies a
FCLT with time-transformed Brownian limit; that is,
%
%
\begin{equation}
\label{arriveFCLT}\quad \hat{N}_n (t) \equiv n^{-1/2}
\bigl(N_n (t) - n \Lambda(t)\bigr) \quad\Rightarrow\quad \hat{N} (t) \equiv
c_{\lambda} B_{\lambda} \bigl(\Lambda(t)\bigr) \qquad\mbox{in }\D
\end{equation}
as $n \ra\infty$, where $B_{\lambda}$ is a standard Brownian motion
(with the subscript $\lambda$ indicating that it is associated with
the arrival process),
$\Lambda(t)$ is the total arrival rate over the interval $[0,t]$, that is,
%
%
\begin{equation}
\label{Lambda} \Lambda(t) \equiv\int_{0}^{t}
\lambda(s) \,ds,
\end{equation}
and $c_{\lambda}^2$ is an arrival-process variability parameter.
A principal case is $N_n$ being a nonhomogeneous Poisson process for
each $n$, in which case $c_{\lambda} = 1$ in (\ref{arriveFCLT}).
Other explicit arrival process models can be constructing from random
or deterministic time-changes of stationary processes (e.g., renewal processes)
known to satisfy a FCLT. For a rate-$1$ renewal process, $c_{\lambda
}^2 = \sigma_{\lambda}^2/m^2_{\lambda} = \sigma_{\lambda}^2$, where
$m_{\lambda} =1$ is the mean and $\sigma_{\lambda}^2$ is the
variance of an interrenewal time; see Section~7.3 of \cite{W02}.

We will specify smoothness assumptions for the model data $(\lambda,s,
G, F)$ in the next section.
These assumptions allow the staffing function $s$ to decrease in OL intervals.
Thus, as discussed in Section~1 of \cite{LW10}, it is important to
consider what happens in the queueing system if the staffing must
decrease when the service facility is full.
Here we simply assume that the required number of customers are forced
out of the system whenever that happens, without having any future
impact on the system,
that is, without altering the queue content or generating subsequent
retrials. Since the service times are exponential, we need not pay attention
to which customers are forced to leave. However, in the next section we
assume that the staffing function is feasible for the fluid model
(which can be achieved
since it is a deterministic system). We say a staffing function is
feasible if no customer is forced out of service (with unfinished
business) when the staffing function decreases. Moreover, we make
conditions ensuring that the staffing function is asymptotically
feasible for the sequence of stochastic models. Hence, any staffing
function infeasibility is asymptotically negligible.

Let $B_n (t, y)$ [$Q_n (t, y)$] denote the number of customers in
service (queue)
at time $t$ that have been so for time at most $y$.
Let $B_n (t) \equiv B_n (t, \infty)$ [$Q_n (t) \equiv Q_n (t, \infty
)$], the total number of customers in service (queue).
Let $X_n (t) \equiv B_n (t) + Q_n (t)$, the total number of customers
in the system.
Let $W_n (t)$ be the head-of-line waiting time (HWT), that is, the
elapsed waiting time for the customer at the head of the line at time $t$
(the customer who has been waiting the longest). Let $V_n (t)$ be the
potential waiting time (PWT) at time $t$, that is,
the virtual waiting time at time~$t$ (the waiting time if there were
a new arrival at time $t$)
assuming that customer never would abandon (but without actually
altering any arrival's abandonment behavior).
Let $A_n (t)$ be the number of abandonments, and let $D_n (t)$ be the
number of departures (service completions) in
the interval $[0,t]$. We can exploit flow conservation to write
%
%
\begin{equation}
\label{aban} A_n (t) = X_n (0) + N_n (t) -
D_n (t) - X_n (t),\qquad t \ge0.
\end{equation}

Let the associated FWLLN-scaled or fluid-scaled processes be
%
%
\begin{eqnarray}
\label{fluidScaled} \bar{B}_n (t, y) & \equiv& n^{-1}
B_n (t, y),\qquad \bar{Q}_n (t, y) \equiv n^{-1}
Q_n (t, y),
\nonumber
\\
\bar{X}_n (t) & \equiv& n^{-1} X_n (t),\qquad
\bar{D}_n (t) \equiv n^{-1} D_n (t),
\\
\bar{A}_n (t) & \equiv& n^{-1} A_n (t),\qquad t
\ge0.\nonumber
\end{eqnarray}
The waiting times $W_n (t)$ and $V_n (t)$ are not scaled in the fluid limit.
Let the associated FCLT-scaled processes be
%
%
\begin{eqnarray}
\label{fcltScaled}\qquad \hat{B}_n (t, y) & \equiv& n^{-1/2}
\bigl(B_n (t, y) - n B (t,y)\bigr),
\nonumber
\\
\hat{Q}_n (t, y) & \equiv& n^{-1/2} \bigl(Q_n
(t, y) - n Q (t, y)\bigr),
\nonumber
\\
\hat{X}_n (t) & \equiv& n^{-1/2} \bigl(X_n (t)
- n X (t)\bigr),\qquad \hat {D}_n (t) \equiv n^{-1/2}
\bigl(D_n (t) - n D (t)\bigr),
\\
\hat{A}_n (t) & \equiv& n^{-1/2} \bigl(A_n (t)
- n A (t)\bigr),\qquad \hat {W}_n (t, y) \equiv n^{1/2}
\bigl(W_n (t) - w (t)\bigr),
\nonumber
\\
\hat{V}_n (t) & \equiv& n^{1/2} \bigl(V_n (t)
- v (t)\bigr),\nonumber
\end{eqnarray}
where $(B (t,y), Q (t,y), X (t), A(t), D(t), w (t), v(t))$ is the
vector of fluid model performance functions,
which will arise as the deterministic limit functions for the associated
FWLLN-scaled processes, already identified in \cite{LW10}.

Our objective is to (i) show that the FWLLN-scaled processes in (\ref
{fluidScaled}) converge in distribution to the previously studied
deterministic fluid model quantities, (ii)~show that the associated
FCLT-scaled processes in (\ref{fcltScaled}) converge in distribution to
a nonstationary zero-mean Gaussian process and (iii) characterize the
dynamics of this Gaussian process
and identify its time-varying variance functions.

\section{The associated deterministic fluid model}\label{secFluid}

The associated deterministic $G_t/M/s_t+\GI$ fluid model depends on the
same model data as the $G_t/M/s_t+\GI$ queueing model
except for the arrival process. The fluid model depends on the arrival
process only through the arrival-rate function $\lambda$.
Thus the fluid model neither captures the full distribution of the
arrival processes
nor the Brownian limit in (\ref{arriveFCLT}). [However, the limit in
(\ref{arriveFCLT}) does affect the FCLT.]
The remaining
functions $(\lambda, s, G, F)$ specify an associated $G_t/M/s_t+\GI$
fluid model as studied in \cite{LW10,LW11a}. All components play an
important role in its performance
description, including the c.d.f. $F$ beyond its mean.

For the fluid model,
$G(x)$ is the proportion of any quantity of fluid that
completes service
by time $x$ after it enters service, and $F(x)$ is the proportion of
any quantity of fluid that abandons
by time $x$ after it enters the queue if it has not already entered service.
We assume that the assumptions for the fluid model
in \cite{LW10} are satisfied here. In \cite{LW10} we uniquely
characterize all the fluid performance functions
under the stated assumptions. We exploit that characterization here, so
that explains how all these assumptions are used.
We conjecture that the limits here can be extended by weakening the conditions,
but we anticipate that will lead to more general, but less tractable,
limits in the FWLLN and FCLT, such as measure-valued functions and
stochastic processes,
as in \cite{KR10,KR11}.

Of special note is the smoothness assumption from \cite{LW10}: we assume
that the functions
$\Lambda$, $s$ and $F$ introduced above are differentiable with
derivatives $\lambda$, $\dot{s}$ and $f$ that are in the space $\C_{\mathrm{pc}}$,
the subspace of $\D$ containing piecewise-continuous functions,
having only finitely many discontinuities in each bounded interval. For
the FCLT in OL intervals, Theorem \ref{thFCLT},
we also assume that $\lambda$ is differentiable as well.
In addition, we assume that
$G^c (x) \equiv1 - G(x) = e^{-\mu x}$, $F^c (x) \equiv1 - F(x) > 0$
for all $x$,
$\lambda_{\mathrm{inf}} \equiv\inf_{0 \le u \le t}{\lambda(u)} > 0$
and $s_{\mathrm{inf}} \equiv\inf_{0 \le u \le t}{s (u)} > 0$.

Consistently with \cite{LW10}, but contrary the terminology for fluid
scaled processes in (\ref{fluidScaled}),
we will denote the fluid performance measures without a bar; thus
$B(t,y)$ [$Q(t,y)$]
denotes the fluid content in service (queue)
at time $t$ that has been so for time at most $y$. These quantities
have densities, that is,
%
%
\begin{equation}
\label{densities} B(t,y) = \int_{0}^{y} b(t,x) \,dx \quad\mbox{and}\quad Q(t,y) = \int_{0}^{y} q(t,x) \,dx.
\end{equation}
Since we have exponential service here, it suffices to focus on the
total fluid content in service
$B(t) \equiv B(t,\infty)$. Let $Q(t) \equiv Q (t, \infty)$ and $X(t)
\equiv B (t) + Q (t)$. Let $w (t)$ be the
head-of-line waiting time (HWT), called the boundary waiting time in
\cite{LW10}; let $v (t)$ be the potential waiting time (PWT) of new
fluid input at time $t$,
both defined essentially the same as $W_n(t)$ and $V_n (t)$ in the
queueing model.

We assume that fluid model starts out underloaded with initial fluid
content $X (0) = B(0)$, where necessarily $B(0) \le s(0)$ and $Q(0) = 0$.
Since the service-time c.d.f. $G$ is exponential, we make no assumption
about the length of time that initial fluid
has been in service.
We assume that the fluid model
has only finitely many switches between underloaded (UL) and overloaded
(OL) intervals in any bounded time interval;
conditions for that property to hold are given in \cite{LW11a}.

The OL and UL intervals are carefully defined in \cite{LW10} (to which
we refer for details).
In this paper, we impose the stronger assumption that the fluid model
is never critically loaded except at the finitely many switching points
in any bounded time interval.
In particular, if $[\tau_1, \tau_2]$ is a UL interval with switching
times at its endpoints, so that $X(\tau_i) = s(\tau_i)$ for $i = 1,2$,
then we require that $X(t) < s(t)$ for all $t$, $\tau_1 < t < \tau_2$.
On the other hand, if $[\tau_1, \tau_2]$ is a OL interval with
switching times at its endpoints,
then we require that $X(t) > s(t)$ for all $t$, $\tau_1 < t < \tau_2$.

The UL intervals are relatively elementary because then the fluid model
is equivalent to an associated infinite-capacity model.
However, the OL intervals are more complicated.
First, as in \cite{LW10}, it is important to assume that the fluid
staffing functions $s$ is feasible, that is, that
its decreasing never forced fluid out of service.
In \cite{LW10} we also show how to construct the minimum feasible
staffing function greater than or equal to any given staffing function.

Here we assume that the flow rate of fluid into service is strictly
positive throughout the OL interval $[\tau_1, \tau_2]$;
that is, we assume that the rate fluid enters service due to new
service capacity becoming available satisfies
%
%
\begin{equation}
\label{flowInto} b(t,0) = s(t)\mu+ \dot{s} (t) \ge b_{\mathrm{inf}} >
0,\qquad
\tau_1 \le t \le\tau_2.
\end{equation}
Together with the FWLLN, condition (\ref{flowInto}) implies that the
probability the staffing function $s_n (t)$ is feasible for the
stochastic model throughout the interval $[\tau_1, \tau_2]$
converges to $1$ as $n \ra\infty$.

We now review the fluid performance functions during an OL interval.
From Section~6 of \cite{LW10}, we know that with $\GI$ service the
fluid density in an overloaded interval requires solving a fixed point
equation, but with $M$ service the service content density during an OL
interval is given explicitly by
%
%
\begin{equation}
\label{fluidDensityB} b(t,x) = b(t-x, 0) G^c (x) 1_{\{x \le t\}} +
b(0,x-t) \frac{G^c (x)}{G^c (x - t)} 1_{\{x > t\}},
\end{equation}
where $G^c (x) \equiv1 - G(x) \equiv e^{-\mu x}$, $b(t, 0) = \dot{s}
(t) + s(t) \mu$, the rate fluid enters service at time $t$,
and $b(0,x)$ is the initial service content density,
part of the initial data.

In \cite{LW10} the queue during an overloaded interval is analyzed by
focusing on the fluid content density $\tilde{q} (t,x)$
assuming no flow into service. Paralleling (\ref{fluidDensityB}),
assuming an initially empty queue, it can be written explicitly as
%
%
\begin{equation}
\label{fluidDensityQno} \tilde{q}(t,x) = \lambda(t-x) F^c (x)
\end{equation}
for $x \le t$, which is all we consider. By Corollary 2 of \cite
{LW10}, the queue content density itself is
%
%
\begin{equation}
\label{fluidDensityQ} q(t,x) = \tilde{q}(t,x)1_{\{x \le w(t) \}} = \lambda(t-x)
F^c (x) 1_{\{x \le w(t)\}},
\end{equation}
so that $q (t,x)$ is simply $\tilde{q} (t, x)$ truncated in the second
variable at its right boundary,
the HWT $w(t)$.

By Theorem 3 of \cite{LW10},
the fluid HWT $w$ is the unique solution to the ODE
%
%
\begin{equation}
\label{wODE} \dot{w} (t) \equiv\frac{d w }{dt} (t) = 1 - \frac{b(t,0)}{\tilde
{q} (t, w(t))} = 1
- \frac{\dot{s} (t) + s(t) \mu}{\lambda
(t-w(t))F^c (w(t))},
\end{equation}
where $b(t, 0) = \dot{s} (t) + s(t) \mu$ is the rate that fluid
enters service.
Our assumptions imply that both the numerator and the denominator in
the fraction in (\ref{wODE})
are strictly positive; thus $-\infty< \dot{w} (t) < 1$ for all $t$ in
the OL interval.
The ODE in (\ref{wODE}) is
equivalent to the integral equation
%
%
\begin{equation}
\label{wIntEq} w(t) = \int_{0}^{t} \biggl(1 -
\frac{b(u,0)}{\tilde{q}(u,
w(u))} \biggr) \,du,\qquad t \ge0.
\end{equation}

By Theorem 5 of \cite{LW10},
the fluid PWT $v(t)$ is as the unique solution of the equation
%
%
\begin{equation}
\label{vEq} v\bigl(t-w(t)\bigr) = w(t) \quad\mbox{or, equivalently}\quad v(t) = w\bigl(t +
v(t)\bigr),
\end{equation}
which can be solved given the BWT $w$. Because of assumption (\ref
{flowInto}), $v$ is
a continuous function. Indeed, both $w$ and $v$ are differentiable
except at only finitely many points.
From (\ref{vEq}), we see that the derivatives are related by
%
%
\begin{equation}
\label{vDeriv}\quad \dot{v} \bigl(t - w(t)\bigr) = \frac{\dot{w} (t)}{1 - \dot{w}(t)}
\quad\mbox{or, equivalently}\quad \dot{v} (t) = \frac{\dot{w} (t + v(t))}{1
- \dot{w} (t + v (t))},
\end{equation}
which is bounded because of condition (\ref{flowInto}).

Since the service is exponential and the service facility is full in an
OL interval, the total fluid departure (service completion) in $[0,t]$
is $D(t) = S(t) \mu$, where $S(t) \equiv\int_{0}^{t} s(u) \,du$.
Finally during an OL interval, the fluid abandonment over $[0,t]$ is
%
%
\begin{equation}
\label{abanFluid} A(t) = \int_{0}^{t} \alpha(s)
\,ds\qquad \mbox{where } \alpha(s) = \int_{0}^{\infty}
Q(s,x) h_F (x) \,dx
\end{equation}
with $h_F (x) \equiv f(x)/F^c (x)$, the hazard rate function associated
with the c.d.f. $F$, which is finite for all $x$
because $f$ is an element of $\D$ and $F^c (x) > 0$ for all $x$.

\section{Heavy-traffic limits during an overloaded interval}\label{secOL}

Recall that the system is said to be OL or UL if the associated fluid
model is OL or UL,
which depends on the model parameters.
The definitions were given in Section~\ref{secFluid}.
We establish the many-server heavy-traffic limits over successive UL
and OL intervals, using the limit at the right endpoint of the
previous interval to provide the limit for the initial conditions
needed in the successive interval, for example, as in \cite{KW90}.
As indicated in the last section, we assume that the fluid model is
initially underloaded. Thus there are UL intervals $[\tau_{2i}, \tau
_{2i + 1}]$, $i \ge0$,
and OL intervals $[\tau_{2i+1}, \tau_{2i + 2}]$, $i \ge0$, with some
finite number of these covering some overall finite time interval of
interest $[0,T]$.
We consider these intervals recursively, referring to each interval in
question as $[0, \tau]$. It should be shifted to the appropriate time.

For the first UL interval, we assume that we have a limit for the
initial conditions, in particular,
%
%
\begin{equation}
\label{InitialLim} \bar{X}_n (0) \Rightarrow X (0) \quad\mbox{and}\quad
\hat{X}_n(0) \Rightarrow\hat {X} (0) \qquad\mbox{in }\RR\mbox{ as } n \ra\infty,
\end{equation}
where $X(0)$ is deterministic with $X (0) \le s(0)$.
For all subsequent intervals, UL and OL, the limit in (\ref{InitialLim})
will hold with $X(0) = s(0)$
as a consequence of the limit in the previous subinterval.

We first consider the more challenging case of an overloaded interval
$[0,\tau]$, assuming limits for the
initial values as in (\ref{InitialLim}), with $X(0) = s(0)$.
We first state the FWLLN. The proofs are given afterward in later sections.
Unlike \cite{LW10}, here we have assumed that the rate of the flow
into service $b(t,0) = s(t) \mu+ \dot{s} (t) > b_{\mathrm{inf}} (\tau) > 0$,
so that the fluid PWT $v$ satisfying (\ref{vEq}) is continuous. Let
$\D
^k$ denote the $k$-fold product space of $\D$ with the associated
product topology.
%
\begin{theorem}[(FWLLN for each OL interval)]\label{thFWLLN}
Consider an OL interval $[0,\tau]$ with no critical loading
except at the endpoints. Suppose that (\ref{InitialLim}) holds
with $X(0) = s(0)$. Then
%
%
\begin{equation}
\label{fLim} (\bar{N}_n, \bar{D}_n,
\bar{X}_n, \bar{Q}_n, \bar{B}_n, \bar
{A}_n, W_n, V_n) \Rightarrow(\Lambda, D,
X, Q, B, A, w, v)
\end{equation}
in $\D^8 ([0, \tau])$ as $n \ra\infty$, where the converging
processes are defined in Section~\ref{secSeq},
the limit $(\Lambda, D, X, Q, B, A, w, v)$ is the vector of continuous
deterministic fluid-model functions defined in Section~\ref{secFluid}
and characterized in \cite{LW10}, having $Q \ge0$, $X = Q + s$
and $B = s$.
\end{theorem}

We next state the associated FCLT establishing the Gaussian refinement
to the fluid approximation in an OL interval. As indicated in the
introduction, we primarily focus on the number in system, $X_n (t)$. We
express the limit for $X_n (t)$ with the general initial conditions in
(\ref{InitialLim}) in terms of an associated limit for the special case
in which all servers are busy, and the queue is empty. Let $X_n^{*}(t)$
be the number in system for the special initial condition in which all
servers are busy and the queue is empty at time $0$, that is,
$X_{n}^{*}(0) = s_{n}(0) = \lceil n\cdot s(0)\rceil$. Let the other
processes associated with this special initial condition be defined
similarly. We now assume that the arrival rate function $\lambda$ is
differentiable in order to work with the partial derivative
%
%
\begin{equation}
\label{partialQ} \tilde{q}_x (t, x) \equiv\frac
{\partial
\tilde{q}(t,x)}{\partial x}.
\end{equation}
Let $\sB$ denote a standard (drift
$0$, diffusion coefficient $1$) Brownian motion (BM). [Recall that
$B(t)$ is already used to denote the fluid content in service.] Let $e$
denote the identify function in $\D$, that is, $e(t) = t$.

\begin{theorem}[(FCLT for each OL interval)]\label{thFCLT}
Consider an OL interval $[0,\tau]$
with no critical loading
except at the endpoints. Assume that the arrival rate function $\lambda
$ is differentiable and the patience p.d.f. $f$ is continuous. Suppose
that (\ref{InitialLim}) holds with $X(0) = s(0)$. Then
%
%
\begin{eqnarray}
\label{gXlim} && \bigl(\hat{N}^{*}_n,
\hat{D}^{*}_n, \hat{X}_n^{*},
\hat{Q}_n^{*}, \hat{B}_n^{*},
\hat{W}^{*}_n, \hat{V}^{*}_n,
\hat{A}^{*}_n, \hat {X}_n\bigr)
\nonumber\\[-8pt]\\[-8pt]
&&\qquad \Rightarrow\bigl(\hat{N}^{*}, \hat{D}^{*},
\hat{X}^{*}, \hat {X}^{*}, 0e, \hat{W}^{*},
\hat{V}^{*}, \hat{A}^{*},\hat{X}\bigr) \qquad\mbox{in }\D^9
\bigl([0,\tau]\bigr),\nonumber
\end{eqnarray}
where the superscript $*$ denotes the special initial condition with
all servers busy and an empty queue,
the converging processes are defined in Section~\ref{secSeq},
and the limit process with the special initial condition, $(\hat
{N}^{*}, \hat{D}^{*}, \hat{X}^{*}, \hat{X}^{*},\hat{W}^{*},\break  \hat
{V}^{*}, \hat{A}^{*})$,
is a mean-zero Gaussian process
having continuous sample paths. If $\hat{X} (0)$ is Gaussian with mean
$0$, then $\hat{X}$ is a mean-zero Gaussian process too.
The limit processes are $\hat{N}^{*} (t) \equiv c_{\lambda} \sB
_{\lambda} (\Lambda(t))$ and $\hat{D}^{*} (t) \equiv\sB_s (D(t))$,
while
%
%
\begin{eqnarray}
\label{Glims} \hat{X} (t) & \equiv& \hat{X}^{*} (t) + \hat{X}(0)
F_w^c (t),\qquad \hat{X}^{*} (t) \equiv\sum
_{i= 1}^{3} \hat{X}^{*}_{i}
(t),
\nonumber
\\
F^c_w (t) & \equiv& e^{- \int_{0}^{t} h_F (w(u)) \,du},\qquad \hat
{X}^{*}_{i} (t) \equiv\int_{0}^{t}
K_{i} (t,u) \,d \sB_{i} (u),
\nonumber\\
\hat{W}^{*} (t) & \equiv& \sum_{i= 1}^{3}
\hat{W}^{*}_{i} (t),\qquad \hat{W}^{*}_{i}
(t) \equiv\int_{0}^{t} J_{i} (t,u) \,d
\sB_{i} (u),
\\
\hat{V}^{*} (t) & \equiv& \frac{\hat{W}^{*} (t + v(t))}{1 - \dot
{w}(t + v(t))},
\nonumber
\\
\hat{A}^{*} (t) & \equiv& \hat{N}^{*} (t) -
\hat{D}^{*} (t) - \hat {X}^{*} (t),\qquad t \ge0,
\nonumber
\end{eqnarray}
where $h_F (x) \equiv f(x)/F^c (x)$ is the patience hazard rate, $w(t)$
is the fluid HWT, $v(t)$ is the fluid PWT,
$\sB_1 \equiv\sB_{\lambda}$, $\sB_2 \equiv\sB_{s}$ and $\sB_3
\equiv\sB_{a}$ are independent (standard) BMs,
%
%
\begin{eqnarray}
\label{kernelsH} H(t,u) & \equiv& \exp{ \biggl\{\int_{u}^{t}
h(v) \,dv \biggr\}},
\nonumber
\\
h(t) & \equiv& \frac{b(t,0) \tilde{q}_x (t, w(t))}{\tilde{q}^2 (t, w(t))} =\bigl(1 - \dot{w} (t)\bigr)
\frac{ \tilde{q}_x (t, w(t))}{\tilde{q} (t,
w(t))}
\\
& = & \bigl(1 - \dot{w} (t)\bigr) \biggl(\frac{- \dot{\lambda} (t -
w(t))}{\lambda(t - w(t))} - h_{F}
\bigl(w(t)\bigr) \biggr)\nonumber
\end{eqnarray}
and
%
%
\begin{eqnarray}
\label{kernels}
J_i (t, u) & \equiv& I_i (u) H(t,u),
\nonumber
\\
I_{1} (u) & \equiv& I_{\lambda} (u) \equiv\frac{ c_{\lambda} \sqrt {F^c
(w(u))b(u,0)}}{\tilde{q}(u, w(u))},\nonumber\\
\bar{I}_{1} (u) &\equiv&\frac{
c_{\lambda} F^c (w(u))b(u,0)}{\tilde{q}(u, w(u))},
\nonumber\\[-8pt]\\[-8pt]
K_{1} (t, u) & \equiv& K_{\lambda} (t, u) \nonumber\\
&\equiv&
c_{\lambda} F^c (t-u) \sqrt{\lambda(u)} {\mathbf1}_{\{t- w(t) \le u \le t\}}
\nonumber
\\
&&{} + \tilde{q}\bigl(t, w(t)\bigr) \sqrt{\lambda(u)} \bar {I}_{1}
\bigl(L^{-1}(u) \bigr) H \bigl( t, L^{-1}(u) \bigr) {\mathbf
1}_{\{0 \le u \le t- w(t)\}},
\nonumber
\\
\label{kernels2} I_{2} (u) & \equiv& I_{s} (u) \equiv-
\frac{ \sqrt{b(u,0) - \dot{s}
(u)} }{ \tilde{q}(u, w(u))},
\nonumber
\\
K_{2} (t, u) & \equiv& K_{s} (t, u) \equiv\tilde{q}
\bigl(t, w(t)\bigr) J_s (t, u) \nonumber\\
&=& -\sqrt{b(t,0) - \dot{s}(t)} H(t,u),
\nonumber
\\
I_{3} (u) & \equiv& I_{a} (u) \equiv- \frac{ \sqrt {F(w(u))b(u,0)}}{\tilde{q}(u,
w(u))},\nonumber\\[-8pt]\\[-8pt]
\bar{I}_{3} (u) &\equiv& - \frac{\sqrt{F^c (w(u))F(w(u))}}{\tilde{q}(u, w(u))},
\nonumber\\
\qquad K_{3} (t, u) & \equiv& K_{a} (t, u) \nonumber\\
&\equiv&-\sqrt{
\lambda(u) F(t-u) F^c (t-u)} {\mathbf1}_{\{t- w(t) \le u \le t\}}
\nonumber
\\
&&{} + \tilde{q}\bigl(t, w(t)\bigr) \sqrt{\lambda(u)}\bar {I}_{3}
\bigl(L^{-1}(u) \bigr)H \bigl( t, L^{-1}(u) \bigr) {\mathbf
1}_{\{0 \le u \le t- w(t)\}}
\nonumber
\end{eqnarray}
and $L^{-1}$ is the inverse of the function $L(t) = t-w(t)$.
The limit process $\hat{W}^{*}$ is also characterized
as the unique solution to the SDE
%
%
\begin{equation}
\label{sde} d \hat{W}^{*} (t) = h(t) \hat{W}^{*} (t) \,dt
+ I(t) \,d \sB(t)
\end{equation}
for $\sB$ a BM, $h(t)$ in (\ref{kernels}) and
%
%
\begin{eqnarray}
\label{Isquared} I(t)^2 &\equiv&\sum_{i=1}^{3}
I^2_i (t) \nonumber\\[-9pt]\\[-9pt]
&=& \frac{b(t,0) - \dot{s} (t)
+ [F(w(t)) + c^2_{\lambda} F^c (w(t))] b(t,0)}{\tilde{q}^2(t,
w(t))}.\nonumber
\end{eqnarray}
\end{theorem}

\begin{remark}[(Additivity of variability)]\label{rmAdd}
It is significant that the three sources of randomness appear
additively (independently)
in the limit process $(\hat{X}^{*}, \hat{W}^{*})$ in~(\ref{Glims}).
The arrival process
variability is captured by $(\hat{X}^{*}_1, \hat{W}^{*}_1)$ and the
BM $\sB_1 \equiv\sB_{\lambda}$;
the service-time variability is captured by $(\hat{X}^{*}_2, \hat
{W}^{*}_2)$ and the BM $\sB_2 \equiv\sB_{s}$; while the
patience-time variability is captured by $(\hat{X}^{*}_3, \hat
{W}^{*}_3)$ and the BM $\sB_3 \equiv\sB_a$, where the three BMs
are mutually independent. Moreover, the four separate sources of
randomness, including $\hat{X} (0)$ for the initial condition in (\ref
{InitialLim}),
which is independent of $(\sB_1, \sB_2, \sB_3)$,
appear additively in the limit process $\hat{X}$.

This nice separation of the components of the variability can be
understood by considering the two-parameter process $Q_n (t,y)$, which depicts
the number of customers in the queue at time $t$ with elapsed patience
time at most $y$ in model $n$ during an OL interval.
The arrivals influence this process at $y= 0$, the lower limit of $y$,
because new arrivals have elapsed patience time $0$.
Because of the FCFS service discipline, the flow into service occurs
from the upper limit of $y$, at $y = W_n (t)$; the customers enter from
the head of the queue; that is,
those who have waited the longest enter first. Finally,
the abandonment influences the process throughout the entire region and
is thus not primarily determined by the behavior at the extreme
endpoints. In particular,
the abandonment rate for a customer with elapsed patience time $x$ is
precisely the patience hazard rate $h_F (x) \equiv f(x)/F^{c}(x)$,
which operates at time $t$
for all $x$ satisfying $0 < x < W_n (t)$ and thus $0 < x < w(t)$ in the
fluid limit.
\end{remark}

Except for the process $\hat{X}_n (t)$, representing the scaled number
in system, Theorem \ref{thFCLT}
states conclusions about the various processes for the special initial
condition,
with all servers busy but no queue.
From Theorem \ref{thFCLT}, we can deduce a corresponding FCLT for the
other processes with the general initial condition in (\ref{InitialLim}),
provided that we exclude the interval endpoints. Recall that
convergence to a continuous limit in $\D^k ((0, \tau))$ is equivalent to
uniform convergence over each compact subinterval $[t_1, t_2]$ with
$0\!<\!t_1\!<\!t_2\!<\!\tau$.

\begin{theorem}[{[Limits for other processes under
(\ref{InitialLim})]}]\label{thOther}
Under the assumptions of Theorem \ref{thFCLT}, all the processes
with the initial conditions
in (\ref{InitialLim}) converge in the space $\D((0,\tau))$; in
particular,
%
%
\begin{equation}
\label{gXlims2}\qquad (\hat{X}_n, \hat{Q}_n,
\hat{B}_n, \hat{W}_n, \hat{V}_n, \hat
{A}_n) \Rightarrow(\hat{X}, \hat{X}, 0e, \hat{W}, \hat{V}, \hat
{A}) \qquad\mbox{in }\D^6 ((0, \tau)),\vadjust{\goodbreak}
\end{equation}
where $\hat{X}$ is given above in (\ref{Glims}),
%
%
\begin{eqnarray}
\label{Glims2} \hat{V} (t) & \equiv& \hat{V}^{*} (t) +
\frac{\hat{X} (0) F^c_w (t
+ v (t))}{s (t + v(t))\mu+ \dot{s} (t+v(t))}
\nonumber
\\
& = &\frac{\tilde{q} (t + v(t), v(t)) \hat{W}^{*} (t + v(t)) +
\hat{X} (0) F^c_w (t + v (t))}{s (t + v(t))\mu+ \dot{s} (t+v(t))},
\nonumber\\[-8pt]\\[-8pt]
\hat{W} (t) &\equiv& \bigl(1 - \dot{w}(t)\bigr)\hat{V} \bigl(t-v(t)\bigr) = \hat
{W}^{*} (t) + \frac{\hat{X} (0) F^c_w (t)}{\tilde{q} (t, w(t))},
\nonumber
\\
\hat{A} (t) & \equiv& \hat{N}^{*}
(t) - \hat{D}^{*} (t) - \hat{X} (t) + \hat{X} (0).\nonumber
\end{eqnarray}
At the interval endpoints $t = 0$ and $t = \tau$, there is the limit
in $\RR^4$
%
%
\begin{equation}
\label{gXlims0} \bigl(\hat{X}_n (t), \hat{Q}_n (t),
\hat{B}_n (t), \hat{V}_n (t) \bigr) \Rightarrow \biggl(
\hat{X} (t), \hat{X} (t)^+, \hat{X} (t)^-, \frac{\hat{X} (t)^+}{s(t)\mu+
\dot{s}(t)} \biggr).\hspace*{-35pt}
\end{equation}
Consequently, for $t$ an interval endpoint, if $P(\hat{X} (t) < 0) > 0$,
then there is no FCLT for $\hat{Q}_n$ and $\hat{V}_n$ in $\D([0,
\tau))$; if $P(\hat{X} (t) > 0) > 0$,
then there is no FCLT for $\hat{B}_n$ in $\D([0, \tau))$.
\end{theorem}

\begin{remark}[(Switching points)]\label{rmDisc}
We get limits like (\ref{gXlims0}) above and (\ref{BzeroLim}) in Theorem
\ref{thUL}
at all switching points. However,\vspace*{1pt} unlike the limit process $\hat{X}
(t)$ for the scaled number in system, $\hat{X}_n (t)$, which has
continuous sample paths,
the resulting limit processes for the other scaled processes $\hat
{Q}_n (t)$, $\hat{B}_n (t)$ and $\hat{V}_n (t)$, obtained by
combining (\ref{gXlims2})
and (\ref{gXlims0}),
will typically have sample paths that are neither left continuous nor
right continuous at the switching points.
In particular, the failure to have convergence at the left endpoint $0$
in (\ref{gXlims2}) occurs because, under the stated condition,
the limit process would need to have a discontinuity point
at the left endpoint, which is not allowed in the space $\D$. If the
switching point occurred at time $\tau$ within a larger interval, then
convergence could be obtained
in the open interval $(0, \infty)$ in the $M_1$ topology, after
redefining the limits at the switching points, but not the $J_1$
topology; see Chapter~12 of \cite{W02}.
In any case, there are limits at the switching points, but the limit
process obtained for each $t$
typically has discontinuities at all switching points. However,\vspace*{2pt} this
difficulty does not occur for the scaled number in system $\hat{X}_n
(t)$; it has a
continuous limit process, as given in Theorem \ref{thFCLT}.
\end{remark}

Practical engineering approximations can be based on the resulting
Gaussian approximations,
for which we need the time-dependent variances, to go with the
time-varying means provided by the fluid limit.
The key process is $\hat{X}$, so we are primarily interested in the
variance $\operatorname{Var}(\hat{X} (t))$, denoted by $\sigma^2_{\hat{X}} (t)$.
Let $\sigma^2_{\hat{X}^{*}} (t) \equiv \operatorname{Var}(\hat{X}^{*} (t))$ and
$\sigma^2_{\hat{X}^{*}, \hat{W}^{*}} (t) \equiv \operatorname{Cov}(\hat{X}^{*}
(t), \hat{W}^{*} (t))$, and similarly for the other processes.\vadjust{\goodbreak}
%
\begin{corollary}[(Variances)]\label{corVar} Consider an OL
interval $[0,\tau]$ satisfying (\ref{InitialLim}).
The variances and covariances are
\begin{eqnarray*}
\sigma^2_{\hat{X}} (t) & = & \sigma^2_{\hat{X}^{*}}
(t) + \operatorname{Var}\bigl(\hat {X} (0)\bigr) \bigl(F^c_w (t)
\bigr)^2,
\\
\sigma^2_{\hat{X}^{*}} (t) & = & \sum
_{i=1}^{3} \sigma^2_{\hat
{X}^{*}_i} (t) =
\int_{0}^{t} \sum_{i=1}^{3}
K_i (t,u)^2 \,du,
\\
& = & \int_{t-w(t)}^{t}\lambda(s)
F^{c}(t-s) \bigl(C_{\lambda}^{2} F^{c}(t-s)+F(t-s)
\bigr)\,ds
\\
&&{} + \tilde{q}^{2}\bigl(t,w(t)\bigr) \sigma^{2}_{\hat
{W}^{*}}(t),
\\
\sigma^2_{\hat{W}^{*}} (t) & = & \sum
_{i=1}^{3} \sigma^2_{\hat
{W}^{*}_i} (t) =
\int_{0}^{t} \sum_{i=1}^{3}
J_i (t,u)^2 \,du = \int_{0}^{t}
H^2 (t,u) I^2 (u) \,du,
\\
\sigma^2_{\hat{V}^{*}} (t) & = &
\frac{\sigma^2_{\hat{W}^{*}} (t +
v(t))}{(1 - \dot{w} (t + v(t)))^2},
\\
\sigma^2_{\hat{V}} (t) & = &
\sigma^2_{\hat{V}^{*}} (t) + \frac
{\operatorname{Var}(\hat{X} (0)) (F^c_w (t + v (t)))^2}{(s (t + v(t))\mu+ \dot{s}
(t+v(t)))^2},
\\
\sigma^2_{\hat{W}} (t) & = & \bigl(1 -
\dot{w}(t)\bigr)^2 \sigma^2_{\hat{V}} \bigl(t-v(t)
\bigr)
\\
& = &\sigma^2_{\hat{W}^{*}} (t) + \frac{\operatorname{Var}(\hat{X} (0))
(F^c_w (t))^2}{\tilde{q} (t, w(t))^2},
\\
\sigma^2_{\hat{X}^{*}, \hat{W}^{*}} (t) & = & \sum
_{i=1}^{3} \sigma ^2_{\hat{X}^{*}_i, \hat{W}^{*}_i} (t)
= \int_{0}^{t} \sum_{i=1}^{3}
J_i (t,u) K_i (t,u) \,du,
\end{eqnarray*}
where $K_i$, $J_i$, $H$ and $I$ are given in (\ref{kernels}) and
$F^c_w$ is given in (\ref{Glims}).
\end{corollary}

\section{Heavy-traffic limits during an underloaded interval}\label{secUL}

We now consider the easier case of the UL intervals.
As before, we assume convergence of the initial values, as in (\ref
{InitialLim}).
Clearly, $\bar{X}_n (0) \ge0$, so that necessarily $X (0) \ge0$.
For the initial interval, we can have any nonnegative deterministic
value for $X (0)$, provided that $X (0) \le s(0)$.
For all subsequent UL intervals, the limit over the previous OL
interval will force $X (0) = s(0)$.

As before, we focus on $X_n$ instead of $B_n$, because after the
initial interval we can have $X_n (0) > s_n(0)$,
whereas we necessarily have $B_n (0) \le s_n (0)$.
The important observation here is that, under our assumption that there
is no critically loading in the fluid model
except at the switching points, in each UL interval the processes $\bar
{X}_n$ and $\hat{X}_n$
are asymptotically equivalent to
the associated processes $\bar{X}^{\infty}_n$ and $\hat{X}^{\infty
}_n$ in the associated $G_t/M/\infty$ infinite-server model
with the same arrival process, service times and initial conditions,
$X_n^{\infty} (0) \equiv X_n (0)$.
Thus we can apply many-server heavy-traffic (MSHT) limits established
for that model in \cite{PW10}; also see \mbox{\cite{B67,KP97,RT09}}.
(Previous references suffice here; the full force of \cite{PW10} is
only needed to treat the more general $G_t/\GI/\infty$
model associated with OL intervals.)

For the infinite-server model, we can separate the new arrivals from
the customers initially in the system at time $0$.
Since there are infinitely many servers,
these customers do not interact when they enter service. Moreover, by
the Brownian limit in FCLT in (\ref{arriveFCLT}), the arrivals after any
time $t$ are
asymptotically independent of the arrivals before that time $t$. To
treat the new arrivals, we can assume that the system starts empty.
We use a subscript $e$ to denote quantities associated with the system
starting \textit{empty}, and we use the subscript $z$ to denote
quantities associated with the initial content at time \textit{zero}.
Let $\| \cdot\|_{a,b}$ denote the uniform norm over the interval
$[a,b]$, with $\| \cdot\|_{b}$ also denoting the case in which $a = 0$.
%
\begin{theorem}[(FWLLN and FCLT for UL interval)]\label{thUL}
Consider a UL interval $[0,\tau]$ under condition (\ref{InitialLim}),
allowing no critical loading except at the interval endpoints. Then
%
%
\begin{equation}
\label{ULlims} \bar{X}_n \Rightarrow X \equiv X_e +
X_z \quad\mbox{and}\quad\hat{X}_n \Rightarrow \hat{X} \equiv
\hat{X}_e + \hat{X}_z \qquad\mbox{in }\D\bigl([0, \tau]\bigr)
\end{equation}
as $n \ra\infty$, where
%
%
\begin{equation}
\label{fluidParts}\quad X_e (t) = \int_{0}^{t}
G^c (t -s) \lambda(s) \,ds \quad\mbox{and}\quad X_z (t) \equiv X (0)
G^c (t),\qquad t \ge0,
\end{equation}
and $\hat{X}_e$ and $\hat{X}_z$ are independent stochastic processes,
with $\hat{X}_e$ being a mean-zero Gaussian diffusion process
satisfying the (SDE)
%
%
\begin{eqnarray}
\label{sdeUL} d \hat{X}_{e} (t) & = & - \mu\hat{X}_{e}
(t) \,dt + c_{\lambda} \,d \sB_{\lambda} \bigl(\Lambda(t)\bigr) - d
\sB_{s} \biggl(\mu\int_0^t
X_e (u) \,d u \biggr)
\nonumber\\[-8pt]\\[-8pt]
& \deq& - \mu\hat{X}_{e} (t) \,dt + \sqrt{c^2_{\lambda}
\lambda (t) + \mu X_e (t)} \,d \sB(t),\nonumber
\end{eqnarray}
where $\sB_{\lambda}$, $\sB_s$ and $\sB$ are independent standard
BMs and $\hat{X} (0) \equiv0$.
The limit process associated with the initial conditions is
%
%
\begin{equation}
\label{ULinitBB} \hat{X}_z (t) \equiv\hat{X} (0) G^c (t)
+ \sqrt{X(0)} \sB^0 \bigl(G^c (t)\bigr),
\end{equation}
where $\sB^0$ is a standard Brownian bridge independent of $\hat{X}
(0)$ and the BMs in (\ref{sdeUL}).
Equivalently, the limit process $\hat{X}$ satisfies the single SDE
%
%
\begin{eqnarray}
\label{sdeULone} d \hat{X} (t) & = & - \mu\hat{X} (t) \,dt + c_{\lambda} \,d
\sB _{\lambda} \bigl(\Lambda(t)\bigr) - d \sB_{s} \biggl(\mu\int
_0^t X (u) \,d u \biggr)
\nonumber\\[-8pt]\\[-8pt]
& \deq& - \mu\hat{X} (t) \,dt + \sqrt{c^2_{\lambda}
\lambda(t) + \mu X (t)}\,d \sB(t),\nonumber
\end{eqnarray}
where $\hat{X} (0)$ is given in (\ref{InitialLim}).\vadjust{\goodbreak}

If $\hat{X}(0)$ is a mean-zero Gaussian random variable, then $\hat
{X}_z$ and $\hat{X}$ are mean-zero Gaussian processes with
$\sigma^2_{X} (t) \equiv \operatorname{Var}(\hat{X} (t)) = \sigma^2_{e} (t) +
\sigma^2_{z} (t)$,
%
%
\begin{eqnarray}
\label{Vare} \sigma^2_{e} (t) & \equiv& \operatorname{Var}\bigl(
\hat{X}_e (t)\bigr)
\nonumber\\[-8pt]\\[-8pt]
&=& \bigl(c^2_{\lambda} - 1\bigr)
\int_{0}^{t} \bigl(G^c (t -s)
\bigr)^2 \lambda(s) \,ds + \int_{0}^{t} G^c (t -s)
\lambda(s) \,ds\nonumber
\end{eqnarray}
and
%
%
\begin{equation}
\label{Vari} \sigma^2_{z} (t) \equiv \operatorname{Var}\bigl(
\hat{X}_z (t)\bigr) = X (0) G (t) G^c (t) + \operatorname{Var}\bigl(
\hat{X} (0)\bigr) \bigl(G^c (t)\bigr)^2.
\end{equation}
In addition, $\|\bar{B}_n - \bar{X}_n \|_{\tau} \Rightarrow0$, so that
%
%
\begin{equation}
\label{BfluidLim} (\bar{X}_n, \bar{B}_n,
\bar{Q}_n) \Rightarrow(X, X, 0 e) \qquad\mbox{in }\D ^3 \bigl([0,
\tau]\bigr) \mbox{ as } n \ra\infty,
\end{equation}
while, restricted to the open interval $(0, \tau)$,
%
%
\begin{equation}
\label{BGaussianLim} (\hat{X}_n, \hat{B}_n,
\hat{Q}_n) \Rightarrow(\hat{X}, \hat{X}, 0e) \qquad\mbox{in }\D^3
((0, \tau)).
\end{equation}
At the interval endpoints $t = 0$ and $t = \tau$,
%
%
\begin{equation}
\label{BzeroLim} \bigl(\hat{X}_n (t), \hat{B}_n (t),
\hat{Q}_n (t)\bigr) \Rightarrow\bigl(\hat{X} (t), \hat{X} (t)^-,
\hat{X} (t)^+\bigr) \qquad\mbox{in }\RR^3.
\end{equation}
Consequently, the limit process $\hat{X} (t)$ for the scaled number in
system $\hat{X}_n (t)$ has continuous sample paths, whereas the limit
processes for the scaled number in queue and in service, $\hat{Q}_n
(t)$ and $\hat{B}_n (t)$, typically have sample paths
that are neither left continuous nor right continuous at the switching points.
Thus, if $X(0) = s(0)$ and $P(\hat{X} (0) < 0) > 0$,
then there is no FCLT for $\hat{Q}_n$ in $\D([0, \tau))$; if $X(0) =
s(0)$ and $P(\hat{X} (0) > 0) > 0$,
then there is no FCLT for $\hat{B}_n$ in $\D([0, \tau))$.
\end{theorem}

The remainder of this paper is concerned with proving all the stated results.

\section{\texorpdfstring{Proofs of Theorems \protect\ref{thFWLLN} and
\protect\ref{thFCLT} for overloaded intervals}
{Proofs of Theorems 4.1 and 4.2 for overloaded intervals}}\label{secOLproofs}

The proof of Theorems \ref{thFWLLN} and \ref{thFCLT} is rather long,
so we start by giving a brief overview.
As in Theorem \ref{thFCLT}, we focus on the number in system, $X_n
(t)$. To do so, it
is convenient to first consider the number in system during the OL
interval starting with all servers busy and an empty queue.
Hence, we will initially consider the OL interval under this special
initial condition. We will then establish the limit for $X_n (t)$
with general initial conditions in Section~\ref{secInitialPf}. We do
not use the notation with the superscript $*$ until Section~\ref{secInitialPf}.

In Section~\ref{secSpecialInitial} we show that any idleness right
after time $0$ is asymptotically negligible,
implying that the departure process is asymptotically equivalent to
a nonhomogeneous Poisson process with the rate $s(t)\mu$.
In Section~\ref{secQueueNoPf} we state preliminary results for the
queue-length process ignoring all flow into service;
these results follow directly from the infinite-server results in \cite{PW10}.
In Section~\ref{secQrepPf} we establish important representations for
the queue-length process during the OL intervals, allowing flow into service.
In Section~\ref{secQlimPf} we show that many-server heavy-traffic
limits for the queue-length process follow from corresponding
limits for the HOL waiting times.
In Section~\ref{secWrepPf} we establish an important representation
for the HOL waiting times.
In Sections~\ref{secWfwlln} and \ref{fcltPf}, respectively, we
exploit the results above to prove the FWLLN
and the FCLT, still under the special initial condition.
Finally, in Section~\ref{secInitialPf} we prove that corresponding
limits hold for the general initial condition in (\ref{InitialLim}).

\subsection{Arrivals and departures with the special initial
condition}\label{secSpecialInitial}

We start by considering the special initial condition with all servers
busy and an empty queue.
Since we are in an OL interval with $\Lambda(t) > D(t)$ for all $t$,
$0 < t < \tau$,
with the initial net input rate to service $\lambda(0) - s(0) \mu-
\dot{s} (0) > 0$ and the abandonment hazard rate bounded above,
even though some servers could become idle
shortly after time $0$, all
servers become busy and remain busy throughout an interval $[t_{1,n},
t_2]$ for $0 < t_{1,n} = O(1/\sqrt{n}) < t_2 < \tau$.
Thus there are at most $O(\sqrt{n})$ empty servers for a period of
only $O(1/\sqrt{n})$. Thus, the total service completion process
differs from the nonhomogeneous Poisson process with rate $n D(t)$ by
only $O(\sqrt{n}) \times O(1/\sqrt{n}) = O(1)$ as $n \ra\infty$.
Similar reasoning also applies at the right endpoint~$\tau$.
Hence, we can conclude that the
departure (service completion) process satisfies a joint FWLLN with the
arrival process of the form
%
%
\begin{equation}
\label{twoFWLLN} \bigl(\bar{N}_n (t),\bar{D}_n (t)
\bigr) \Rightarrow\bigl(\Lambda(t), D(t)\bigr) \qquad\mbox{in } \D^2\bigl([0,
\tau]\bigr)
\end{equation}
and a corresponding joint FCLT,
%
%
\begin{eqnarray}
\label{twoFCLT} \bigl(\hat{N}_n (t),\hat{D}_n (t)\bigr)
& \Rightarrow& \bigl(\hat{N} (t),\hat{D} (t)\bigr) \qquad\mbox{in } D^2
\bigl([0, \tau]\bigr) \nonumber\\[-8pt]\\[-8pt]
&&\eqntext{\mbox{where }
\hat{N} (t) \equiv c_{\lambda} \sB_{\lambda} \bigl(\Lambda(t)\bigr)
\mbox{ and } \hat{D} (t) \equiv\sB_s \bigl(D(t)\bigr), t \ge0,}
\end{eqnarray}
with $\sB_{\lambda}$ and $\sB_s$ being two independent BMs.

As a consequence of the results above, we determine (relatively
trivial) limits for the number in service, in particular,
%
%
\begin{equation}
\label{lmB} \bar{B}_n \Rightarrow s \quad\mbox{and}\quad\hat{B}_n
\Rightarrow0e \qquad\mbox{in }\D \bigl([0,\tau]\bigr)
\mbox{ as } n \ra\infty.
\end{equation}
As a consequence, we deduce for the number in queue that
%
%
\begin{equation}
\label{XQdiff} \bigl\|\bar{X}_n - (\bar{Q}_n + s)
\bigr\|_{\tau} \Rightarrow0 \quad\mbox{and}\quad\| \hat {X}_n -
\hat{Q}_n\|_{\tau} \Rightarrow0 \qquad\mbox{as } n \ra\infty.
\end{equation}
Hence, to establish limits for $\bar{X}_n$ and $\hat{X}_n$ in $\D
([0,\tau])$, it suffices to focus on $\bar{Q}_n$ and $\hat{Q}_n$,
which is what we do in the following subsections.

\subsection{The queue length ignoring flow into service}\label{secQueueNoPf}

To study the fluid model in overloaded intervals, in \cite{LW10} we
introduced the fluid function $\tilde{q} (t,x)$,
which is the fluid content density in queue, disregarding the flow into
service, that is, under the condition that the flow into service is
turned off.
It is convenient to do the same in order to develop stochastic
refinements. Let $\tilde{Q}_n (t,y)$ be the two-parameter stochastic process
giving the number in queue in model $n$ at time $t$ that have been so
for at most time $y$,
under the condition that the flow into service is turned off. Until
Section~\ref{secInitialPf}, we have the special initial conditions
with all servers busy and an empty queue.

When we turn off all flow into service, the number in service in the
$G_t/M/s_t+\GI$ queueing model with our special initial condition
is asymptotically equivalent to the number in the associated
$G_t/\GI/\infty$ queueing model, starting empty, where the abandonment
c.d.f. $F$ plays the usual role of the service-time c.d.f. in the
infinite-server model.
Hence, we consider the stochastic process $\tilde{Q}_n (t,y)$ in the
queueing model, disregarding flow into service.
Thus, we can apply the FCLT for the
$G_t/\GI/\infty$ queueing model established by \cite{PW10}.

We exploit the representation of $\tilde{Q}_n (t,y)$ from \cite{PW10}.
Let $1_{A} (t)$ be the indicator function of the set $A$, that is,
$1_{A} (t) = 1$ if $t \in A$ and $0$ otherwise.
First, we can write
%
%
\begin{equation}
\label{tildeRep1} \tilde{Q}_n (t,y) = \sum
_{i = N_n ((t-y)-) + 1}^{N_n (t)} {\mathbf 1}\bigl(\tau^n_i
+ \eta_i > t\bigr),\qquad t \ge0, 0 \le y \le t,
\end{equation}
where $\tau^n_i$ is the $i$th arrival time, $\eta_i$ is the $i$th
patience time (the patience time of the arrival at $\tau^n_i$) and
$N_n (t)$ is the arrival counting process in model $n$. The
representation in (\ref{tildeRep1})
is valid because the first $N_n ((t-y)-)$ arrivals will have come
before time $t$.
(The limit process will have continuous sample paths, so that the
consequence of an arrival exactly at time $t$ is asymptotically negligible.)
Hence, the sum in (\ref{tildeRep1}) counts all arrivals in the interval
$[t-y, t]$ who will not have abandoned
by time $t$.

Following \cite{KP97}, the next step in \cite{PW10} is to obtain an
alternative representation
exploiting the sequential empirical process associated with the
successive patience times,
%
%
\begin{equation}
\label{Kn} \bar{K}_n (t, y) \equiv\frac{1}{n} \sum
_{i=1}^{\lfloor nt \rfloor
} {\mathbf1} (\eta_i \le y),\qquad t
\ge0, y \ge0.
\end{equation}
In particular, representation (\ref{tildeRep1}) is equivalent to the
alternative representation
%
%
\begin{equation}
\label{tildeRep2} \tilde{Q}_{n} (t,y) \equiv n \int
_{t-y}^{t} \int_{0}^{\infty}
{\mathbf 1} (x + s > t) \,d \bar{K}_n \bigl(\bar{N}_n (s),
x\bigr)
\end{equation}
for $t \ge0$ and $0 \le y \le t$.
Representation (\ref{tildeRep2}) allows us to exploit the limits $\bar
{K}_n (t,x) \Rightarrow tF(x)$ and
%
%
\begin{equation}
\label{hatKn} \hat{K}_n (t, x) \equiv\sqrt{n}\bigl(
\bar{K}_n (t, x) - F(x)\bigr) \quad\Rightarrow\quad\hat{K} (t, x) \equiv U
\bigl(t, F(x)\bigr)
\end{equation}
in $D([0,\infty), D([0,1], \RR))$, where the limit $\hat{K}$ is a
deterministic transformation of the
standard Kiefer process $U(t,x)$.

From Lemma 2.1 of \cite{PW10}, we obtain the alternative representation
%
%
\begin{eqnarray}
\label{tildeDecomp} \tilde{Q}_n (t,y) & \equiv&\tilde{Q}_{n,1}
(t,y) + \tilde{Q}_{n,2} (t,y) + \tilde{Q}_{n,3} (t,y),
\nonumber
\\
\tilde{Q}_{n,1} (t,y) & \equiv& \sqrt{n} \int_{t - y}^{t}
F^c (t-s) \,d \hat{N}_n (s),
\nonumber\\[-8pt]\\[-8pt]
\tilde{Q}_{n,2} (t,y) & \equiv& \sqrt{n} \int_{t -y}^{t}
\int_{0}^{\infty} {\mathbf1} (x + s > t) \,d
\hat{R}_n (s, x),
\nonumber
\\
\tilde{Q}_{n,3} (t,y) & \equiv& n \int_{t - y}^{t}
F^c (t-s) \lambda (s) \,ds,\nonumber
\end{eqnarray}
where, just as in (2.16) of \cite{PW10},
%
%
\begin{equation}
\label{Rn} \hat{R}_n (t, y) \equiv\sqrt{n} \bar{K}_n
\bigl(\bar{N}_n (t), y\bigr) - \hat{N}_n (t) F(y) -
\sqrt{n} \Lambda(t) F(y)
\end{equation}
with $\bar{K}_n (t, y)$ being the sequential empirical process in
(\ref{Kn}).

Thus, from \cite{PW10} and (\ref{twoFCLT}), it follows that
%
%
\begin{eqnarray}
\label{tildeLim} && \bigl(\hat{Z}_{n,1} (t, y), \hat{Z}_{n,2}
(t, y)\bigr) \equiv n^{-1/2}\bigl(\tilde{Q}_{n,1} (t,y),
\tilde{Q}_{n,2} (t,y)\bigr)
\nonumber\\[-8pt]\\[-8pt]
&&\qquad \Rightarrow\bigl(\hat{Z}_1 (t, y), \hat{Z}_1 (t, y)
\bigr) \qquad\mbox{in }\D ^2\bigl([0,\tau], \D\bigl([0,\infty),
\RR\bigr)\bigr),\nonumber
\end{eqnarray}
jointly with the limit in (\ref{twoFCLT}), where
%
%
\begin{eqnarray}
\label{tildeLimits} \hat{Z}_1 (t, y) & \equiv& \int
_{t - y}^{t} F^c (t-s) \,d \sB
_{\lambda} \bigl(\Lambda(s)\bigr),
\nonumber\\[-8pt]\\[-8pt]
\hat{Z}_2 (t, y) & \equiv& \int_{t -y}^{t}
\int_{0}^{t} {\mathbf1} (x + s > t) \,d \sR(s, x)\nonumber
\end{eqnarray}
with $\sB_{\lambda}$ being a BM and
%
%
\begin{equation}
\label{R} \sR(t, y) \equiv\sK\bigl( \Lambda(t), y\bigr),
\end{equation}
where $\sK(t,y) \equiv\sU(t, F(y))$ and $\sU(t,x)$ is the standard
Kiefer process,
with $(\sK, \sR)$ independent of $\sB_{\lambda}$.
As a consequence, by the continuous mapping theorem with addition,
%
%
\begin{equation}
\label{tildeLimSum}\quad \hat{Z}_n \equiv\hat{Z}_{n,1} +
\hat{Z}_{n,2} \quad\Rightarrow\quad\hat {Z}_1 + \hat{Z}_2
\qquad\mbox{in }\D\bigl([0,\tau], \D\bigl([0,\infty), \RR\bigr)\bigr)
\end{equation}
for $\hat{Z}_i$ in (\ref{tildeLimits}).

From (\ref{tildeLimits}) and (\ref{R}), we see that
the limit process $\hat{Z}_1$ in (\ref{tildeLim}) depends on the
randomness in the arrival process through the BM $\sB_{\lambda}$,
while the limit process $\hat{Z}_2$ in (\ref{tildeLim}) depends
on the randomness in the patience times through $\sR$, and thus the
Kiefer process $\sK$, associated with the abandonment times.
Since flow into service has not yet been considered, the BM $\sB_s$
does not appear yet.
We will exploit this established convergence in (\ref{tildeLim}) in
order to establish our desired FWLLN and FCLT.

\subsection{Representation of the queue-length process}\label{secQrepPf}

We now obtain a representation of the queue-length process $Q_n (t)$ in
this overloaded interval, where now we are allowing
the usual flow into service.
We do so by modifying the representation for $\tilde{Q}_n (t,y)$
constructed above.
In particular, paralleling (\ref{tildeRep1}), for $t > 0$,
we obtain the representation
%
%
\begin{equation}
\label{repQ}\quad\qquad Q_n (t) = \tilde{Q}_n \bigl(t,
W_n (t)\bigr) = \sum_{i = N_n ((t-W_n (t))-) +
1}^{N_n (t)}
{\mathbf1}\bigl(\tau^n_i + \eta_i > t\bigr),\qquad t
> 0.
\end{equation}
We could also obtain a corresponding representation for the
two-parameter process $Q_n (t,y)$, as in (\ref{tildeRep1}), but here we
focus on the
one-parameter processes. The FCFS service discipline is crucial for
obtaining representation (\ref{repQ}); it ensures that
customers enter service from the head of the line.
Representation (\ref{repQ}) does not tell the whole story, however,
because the HOL waiting time $W_n(t)$ remains to be determined.
Moreover, among the first
$N_n ((t - W_n (t))-)$ arrivals, (\ref{repQ}) does not show which
entered service and which abandoned.

Nevertheless, paralleling (\ref{tildeDecomp}) above, we obtain the
alternative representation
%
%
\begin{eqnarray}
\label{qDecomp} Q_n (t) & \equiv& Q_{n,1} (t) +
Q_{n,2} (t) + Q_{n,3} (t),
\nonumber
\\
Q_{n,1} (t) & \equiv& \sqrt{n} \int_{t - W_n (t)}^{t}
F^c (t-s) \,d \hat{N}_n (s),
\nonumber\\[-8pt]\\[-8pt]
Q_{n,2} (t) & \equiv& \sqrt{n} \int_{t -W_n (t)}^{t}
\int_{0}^{\infty} {\mathbf1} (x + s > t) \,d
\hat{R}_n (s, x),
\nonumber
\\
Q_{n,3} (t) & \equiv& n \int_{t - W_n (t)}^{t}
F^c (t-s) \lambda(s) \,ds,\qquad t > 0,\nonumber
\end{eqnarray}
where $\hat{R}_n$ is given in (\ref{Rn}).

\subsection{\texorpdfstring{Limits for $\hat{Q}_n$ given limits for $\hat{W}_n$}
{Limits for Q n given limits for W n}}\label{secQlimPf}

Given limits $W_n \Rightarrow w$ and $\hat{W}_n \Rightarrow\hat{W}$
in $\D([0, \tau])$ for $\hat{W}_n$ in (\ref{fcltScaled}), where $w$
is the
differentiable fluid HWT satisfying the ODE in (\ref{wODE}) and $\hat
{W}$ has continuous sample paths, which we will establish below, we can obtain
limits for $\bar{Q}_n$ and $\hat{Q}_n$ in $\D([0, \tau])$ directly
from the representation in (\ref{qDecomp})
and the limits in (\ref{tildeLim}) by applying the continuous mapping theorem.
In particular,
%
%
\begin{eqnarray}
\label{QlimLLN} \bar{Q}_{n,i} (t) & \equiv& n^{-1}
Q_{n,i} (t) \quad\Rightarrow\quad(0e) (t) \qquad\mbox{for } i = 1,2,
\nonumber\\[-8pt]\\[-8pt]
\bar{Q}_{n,3} (t) & \equiv& n^{-1} Q_{n,3} (t)
\quad\Rightarrow\quad Q_3 (t) \equiv\int_{t-w(t)}^{t}
\lambda(s) F^c(t - s) \,d s\nonumber
\end{eqnarray}
in $\D([0, \tau])$ and
%
%
\begin{eqnarray}
\label{Q1}
&&\hat{Q}_{n,1} (t) \equiv n^{-1/2} Q_{n,1} (t)\nonumber\\
&&\quad\Rightarrow\quad\hat {Q}_1 (t)
\equiv  C_{\lambda}\int_{t - w(t)}^{t}
F^c (t-s) \,d \tilde {\sB}_{\lambda} \bigl(\Lambda(s)\bigr)
\\
&&\hphantom{\quad\Rightarrow\quad\hat {Q}_1 (t)} \equiv
C_{\lambda}\int_{t - w(t)}^{t}
F^c (t-s) \sqrt {\lambda(s)} \,d \sB_{\lambda} (s),
\nonumber
\\
\label{Q2}
&&\hat{Q}_{n,2} (t) \equiv n^{-1/2} Q_{n,2} (t)\nonumber\\
&&\quad\Rightarrow\quad\hat {Q}_2 (t)
\equiv \int_{t -w(t)}^{t} \int_{0}^{t}
{\mathbf1} (x + s > t) \,d \sR(s, x)
\nonumber\\[-8pt]\\[-8pt]
&&\hphantom{\quad\Rightarrow\quad\hat {Q}_2 (t)}
\deq -\int_{t-w(t)}^{t} \sqrt{F(t-s)
F^c (t-s)} \,d \tilde {\sB}_a \bigl(\Lambda(s)\bigr)
\nonumber\\
&&\hphantom{\quad\Rightarrow\quad\hat {Q}_2 (t)}
\deq -\int_{t-w(t)}^{t} \sqrt{F(t-s)
F^c (t-s) \lambda(s)} \,d \sB_a (s),
\nonumber
\\
\label{Qlim}
\qquad&&\hat{Q}_{n,3} (t) \equiv n^{-1/2} \bigl(Q_{n,3}
(t) - n Q(t)\bigr)
\quad
\Rightarrow \quad\hat{Q}_3 (t) \equiv q\bigl(t, w(t)\bigr) \hat{W}
(t),
\nonumber
\\
&&\hat{Q}_{n} (t) \equiv \hat{Q}_{n, 1} (t) +
\hat{Q}_{n, 2} (t) + \hat{Q}_{n, 3} (t)
\\
&&\quad \Rightarrow \quad\hat{Q} (t)\equiv\hat{Q}_1 (t) +
\hat{Q}_2 (t) + \hat{Q}_3 (t) \qquad\mbox{in }\D((0, \tau)),
\nonumber
\end{eqnarray}
where the three limit processes in the last line are independent. This
is not entirely obvious because $\hat{Q}_3$ involves $\hat{W}$, which
in turn involves
the two BMs $\sB_{\lambda}$ and $\sB_a$ appearing in $\hat{Q}_1$
and $\hat{Q}_2$. However, a close observation reveals that $\hat
{Q}_{1}$ and $\hat{Q}_{2}$ involve the two BMs $\tilde{\sB
}_{\lambda}$ and $\tilde{\sB}_{a}$ from time $\Lambda(t-w(t))$ to
time $\Lambda(t)$, according to the representations in (\ref{Q1}) and
(\ref{Q2}); on the other hand, we will see from (\ref{WB}) of Section~\ref{secWcharPf} that $\hat{W}$ involves $\tilde{\sB}_{\lambda}$
and $\tilde{\sB}_{a}$ from time $\Lambda(0)=0$ to time $\Lambda
(t-w(t))$, which thus concludes the independence.
After we establish the limit for $\hat{W}$, we can appropriately group
the terms and separate these three independent BMs. The representation
in Theorem \ref{thFCLT} will thus follow.

We now justify the convergence just stated above. We start with the FWLLN.
The separate FWLLNs for $N_n$, $Z_{n,i}$ and $W_n$ obtained from (\ref
{twoFWLLN}), (\ref{tildeLim}) and by assumption to deterministic limits
imply the joint FWLLN.
Since we divide by $n$, the terms $\bar{Q}_{n,1}$ and $\bar{Q}_{n,2}$
obtained from (\ref{Qlim}) and $\bar{X}_n (0)^{+}$ become
asymptotically negligible.
Using the assumed FWLLN for $W_n (t)$, we can apply the continuous
mapping theorem with the composition map, specifically Theorem 2.4 of
\cite{TW09},
which extends continuity properties for composition maps to the
two-parameter setting, to the second ($y$) coordinate of $\tilde
{Q}_{n,3} (t,y)$
to obtain $\bar{Q}_{n,3} \Rightarrow Q$ in Theorem~\ref{thFWLLN},
which implies that $\bar{Q}_{n} \Rightarrow Q$.

We now turn to the FCLT refinement. Given the FCLT jointly for $N_n$,
$Z_{n,i}$ and $W_n$ obtained from (\ref{twoFCLT}), (\ref{tildeLim}),
again we can apply
the continuous mapping theorem with the composition map in Theorem 2.4
of \cite{TW09}, applied to the second ($y$) coordinate of $\hat
{Z}_{n,i} (t,y)$
in (\ref{tildeLim})
to obtain the desired conclusions for
$\hat{Q}_{n,i} (t)$, for $i = 1,2$.
Note that we only need the FWLLN for $W_n (t)$ for this step; we do not
need the more involved Theorem 2.5 of \cite{TW09}.
From this step, we obtain the convergence of the vector processes, that is,
%
%
\begin{equation}
\label{vecLim}(\hat{N}_n, \hat{D}_n,
\hat{Z}_{n,1}, \hat{Z}_{n,2}, \hat{W}_n,
\hat{Q}_{n,1},\hat{Q}_{n,2}) \quad\Rightarrow\quad (\hat{N}, \hat{D},
\hat{Z}_{1}, \hat{Z}_{2}, \hat{W}, \hat {Q}_{1},
\hat{Q}_{2}).\hspace*{-34pt}
\end{equation}

Next, we treat $\hat{Q}_{n,3}$ in (\ref{Qlim}) by noting
that
%
%
\begin{eqnarray}
\label{q3first}\qquad \hat{Q}_{n,3} (t) & = & \sqrt{n} \biggl(\int
_{t - W_n (t)}^{t} F^c (t-s) \lambda(s) \,ds -
\int_{t - w (t)}^{t} F^c (t-s) \lambda(s) \,ds \biggr)
\nonumber\\[-8pt]\\[-8pt]
& = & \sqrt{n} \biggl(\int_{t - W_n (t)}^{t - w (t)}
F^c (t-s) \lambda (s) \,ds \biggr),\nonumber
\end{eqnarray}
so that we can exploit the continuity of the integrand $\tilde{q} (t,
t - s) = F^c (t-s) \lambda(s)$ to deduce that
%
%
\begin{equation}
\label{diffQ3}\quad \sup_{0 \le t \le\tau}{\bigl\{\bigl|\hat{Q}_{n,3} (t)
- \hat{W}_n (t) \tilde {q}\bigl(t, w(t)\bigr)\bigr|\bigr\}} = o\bigl(\|
\hat{W}_n \|_{\tau}\bigr) \qquad\mbox{as } n \ra\infty,
\end{equation}
so that $\hat{Q}_{n,3} \Rightarrow\hat{Q}_3$ in $\D([0,\tau])$
jointly with the limit in (\ref{vecLim}) for $\hat{Q}_3 (t) \equiv
\tilde{q}(t, w(t)) \hat{W}(t) $
if $\hat{W}_n \Rightarrow\hat{W}$ in $\D([0, \tau])$. Given that
joint convergence,
we can apply the continuous mapping theorem with addition to obtain the limit
$\hat{Q}_{n} \Rightarrow\hat{Q}$ jointly with the other processes,
as stated in the final line of (\ref{Qlim}).

\subsection{\texorpdfstring{Representation of the HOL waiting times $W_n(t)$}
{Representation of the HOL waiting times W n(t)}}
\label{secWrepPf}

It thus remains only to treat the waiting times.
Paralleling the proof of Theorem 3 of \cite{LW10}, we treat the HWT
process $W_n (t)$ by equating two
different expressions for the number of customers to enter service in
an interval $[t, t+\varepsilon]$, where $\ep$ is a small positive number.
Let $E_n (t)$ be the number of customers to enter service in the
interval $[0,t]$. On the one hand, since the fluid model is overloaded
with $\Lambda(t) > D(t)$ for all $t$, $0 < t < \tau$,
the number of customers to enter service is asymptotically equivalent
to the new capacity made available by departures and changes in the
staffing; that is,
as $n \ra\infty$,
%
%
\begin{equation}
\label{enter1} \sup_{0 \le t \le\tau}{\bigl\{\bigl|E_n (t) -
\bigl(D_n (t) + \bigl\lceil n s(t)\bigr\rceil- \bigl\lceil n s (0)
\bigr\rceil\bigr)\bigr|\bigr\}} = o(\sqrt{n}).
\end{equation}
Let $\bar{E}_n (t) \equiv E_n (t)/n$ and $\hat{E}_n (t) \equiv\sqrt {n}(\bar{E}_n (t) - E(t))$ be the associated FWLLN and FCLT scaled processes,
where $E(t) \equiv D(t) + s(t) - s(0)$. It follows from (\ref{enter1})
and the FCLT for $D_n$ in (\ref{twoFCLT}) that
%
%
\begin{equation}
\label{Elim} \bar{E}_n (t) \Rightarrow E (t) \quad\mbox{and}\quad
\hat{E}_n \Rightarrow\hat{E} \qquad\mbox{in }\D\mbox{ as } n \ra\infty,
\end{equation}
where
%
%
\begin{equation}
\label{Elimit} \hat{E} (t) = \hat{D} (t) = \sB_s \bigl(D(t)\bigr),\qquad t
\ge0,
\end{equation}
as in (\ref{twoFCLT}).

On the other hand, the flow into service most come from customers
leaving the queue. Because the
service discipline is FCFS, that flow must come from the customers who
have been in service the longest.
We can again use representation (\ref{repQ}) to represent the flow into
service over an interval.
Let $E_n (t, \varepsilon) \equiv E_n (t+ \varepsilon) - E_n (t)$ and
similarly for the other processes.
As in the proof of Theorem 3 of \cite{LW10},
if we make the interval short enough, then the abandonments will be
asymptotically negligible.
Thus, paralleling equation (28) in \cite{LW10} for the fluid model,
from (\ref{repQ}) we obtain
%
%
\begin{eqnarray}
\label{in1} E_n (t, \varepsilon) & = & I_n (t, \varepsilon) -
A^I_n (t, \varepsilon)
\nonumber\\[-8pt]\\[-8pt]
&&\eqntext{\mbox{where }\displaystyle I_n (t, \varepsilon) \equiv
\sum_{i = N_n ((t - W_n (t))-) +
1}^{N_n(t+ \varepsilon- W_n (t+\varepsilon))}
{\mathbf1}\bigl(\tau^n_i + \eta_i > t\bigr);}
\end{eqnarray}
that is, $I_n (t, \varepsilon)$ is the number of customers removed from
the right boundary of the queue in the time interval $[t, t+\varepsilon]$,
and $A^I_n (t, \varepsilon)$ is the number of those $I_n (t, \varepsilon)$
customers that actually abandon.
Note that $W_n (t+ \varepsilon) \le W_n (t) + \varepsilon$ because the
waiting time of each customer that remains in queue increases at rate $1$.
Hence the upper limit of summation in (\ref{in1}) always is greater than
or equal to the lower limit of summation there.

We now want to show that $A^I_n (t, \ep)$ is appropriately
asymptotically negligible relative to $I_n (t, \ep)$.
For that purpose, observe that
%
%
\begin{equation}
\label{inAban}\quad 0 \le A^I_n (t, \varepsilon) \le
J_n (t, \ep) \equiv\sum_{i = N_n ((t
- W_n (t))-) + 1}^{N_n(t+ \varepsilon- W_n (t+\varepsilon))}
{\mathbf1}\bigl(t < \tau^n_i + \eta_i \le t+
\ep\bigr);
\end{equation}
that is, $J_n (t, \ep)$ is the number of customers in the system at
time $t$, but not at time $t + \ep$,
who would abandon before time $t + \ep$ if they do not enter service first
in the interval $[t, t+\ep]$. The remaining $I_n (t, \ep) - J_n (t,
\ep)$ customers necessarily enter service in the interval $[t, t + \ep]$
because they would not abandon before time $t + \ep$.

We now show that the bound $J_n (t, \ep)$ in (\ref{inAban}) is
asymptotically negligible\vspace*{1pt} relative to $I_n (t, \ep)$ as $\ep
\downarrow0$,
uniformly in $n$ and $t$,
so that we can ignore $A^I_n (t, \ep)$ by choosing $\ep$ suitably small.
We prove that by bounding $J_n (t, \ep)$ above. First, we observe that
$0 \le\tau^n_i \le\tau$ for the arrival times $\tau^n_i$ under
consideration. Thus
%
%
\begin{equation}
\label{Jbd1} J_n (t, \ep) \le I_n (t, \ep) \sup{\bigl
\{P\bigl(t \le\tau^n_i + \eta_i \le t +
\ep| t \le\tau^n_i + \eta_i\bigr)\bigr\}},
\end{equation}
where
%
%
\begin{eqnarray}
\label{Jbd2} && \sup{\bigl\{P\bigl(t \le\tau^n_i +
\eta_i \le t + \ep| t \le\tau^n_i +
\eta_i\bigr)\bigr\}}
\nonumber\\[-8pt]\\[-8pt]
&&\qquad \le\sup_{0 \le t \le\tau}{\bigl\{F^c (t) - F^c
(t + \ep)\bigr\}} \le\|f\|_{\tau} \ep+ o(\ep)
\qquad\mbox{as }\ep\downarrow0,\nonumber
\end{eqnarray}
where $\|f\|_{\tau} < \infty$ because the c.d.f. $F$ has the density
$f$, which has been assumed to be in $\C_{\mathrm{pc}} \subseteq\D$.
To summarize,
%
%
\begin{equation}
\label{Jbd3} A^I_n (t, \ep) \le J^n (t,
\ep) \le K \ep I^n (t, \ep)
\end{equation}
for some constant $K$ (depending on the c.d.f. $F$ and $\tau$) for all
$\ep$ suitably small, uniformly in $n$ and $t$.

We can characterize the asymptotic behavior of the HWT process $W_n
(t)$ by equating the two expressions for $E_n (t, \ep)$ from
(\ref{enter1}) and (\ref{in1}).
Here we act as if the system is always overloaded, and thus use the
infinite-server model representation;
as in (\ref{enter1}), the error in this step is asymptotically negligible.
Now, reasoning as in (\ref{tildeRep1})--(\ref{tildeDecomp}), we
obtain an
alternative representation for $I_n (t, \varepsilon)$ in (\ref{in1}).
In particular,
%
%
\begin{equation}
\label{in2} I_{n}(t,\varepsilon) = n\int_{t-W_{n}(t)}^{t+\varepsilon-W_{n}(t+\varepsilon
)}
\int_{0}^{\infty}{\mathbf1}(s+x>t)\,d\bar{K}_{n}
\bigl(\bar{N}_{n}(s),x\bigr),
\end{equation}
where $\bar{K}_{n} (t,x)$ again is the sequential empirical process in
(\ref{Kn}), and then
%
%
\begin{equation}
\label{in3} I_{n}(t,\varepsilon) = I_{n,1}(t,\varepsilon) +
I_{n,2}(t,\varepsilon) + I_{n,3}(t,\varepsilon),
\end{equation}
where
%
%
\begin{eqnarray}
\label{Idecomp} I_{n,1}(t,\varepsilon) & = & \sqrt{n}\int
_{t-W_{n}(t)}^{t+\varepsilon
-W_{n}(t+\varepsilon)}F^c(t-s)\,d
\hat{N}_{n}(s),
\nonumber
\\
I_{n,2}(t,\varepsilon) & = & \sqrt{n}\int_{t-W_{n}(t)}^{t+\varepsilon
-W_{n}(t+\varepsilon)}
\int_{0}^{\infty}{\mathbf1}(s+x>t)\,d\hat
{R}_{n}(s,x),
\\
I_{n,3}(t,\varepsilon) & = & n \int_{t-W_{n}(t)}^{t+\varepsilon
-W_{n}(t+\varepsilon)}F^c(t-s)
\lambda(s)\,ds,\nonumber
\end{eqnarray}
where $\hat{R}_n$ is from (\ref{Rn}).

\subsection{\texorpdfstring{Proof of Theorem \protect\ref{thFWLLN}: The FWLLN}
{Proof of Theorem 4.1: The FWLLN}}
\label{secWfwlln}

We now prove the FWLLN, still under our special\vspace*{1pt} initial conditions
imposed in Section~\ref{secSpecialInitial}.
We have\vspace*{1pt} $(\bar{N}_n,\bar{D}_n, \bar{Z}_n) \Rightarrow(\Lambda, D,
\tilde{Q})$ in $D([0,\tau])^2 \times D([0,\tau], D([0, 1], \RR))$
for $\bar{N}_n (t) \equiv n^{-1} N_n (t)$ and $\bar{D}_n (t)\equiv
n^{-1} D_n (t)$ in (\ref{fluidScaled}) and $\bar{Z}_n (t,y) \equiv
n^{-1} \tilde{Q}_n (t,y)$
in (\ref{tildeRep1})--(\ref{tildeLim}),
where $(\Lambda, D, \tilde{Q})$ are the components of the fluid model
in Section~\ref{secFluid},
based on the FCLTs in (\ref{arriveFCLT}), (\ref{twoFCLT}) and (\ref
{tildeLim}).
As shown above, we also obtain the FWLLN for $\bar{Q}_n$ once we
obtain the FWLLN for $W_n$.

We now prove the FWLLN for $W_n$; that is, $W_n \Rightarrow w$.
We prove the FWLLN for $W_n$ by applying the compactness approach, as
in Section~11.6 of \cite{W02}.
In particular, we show that the sequence $\{W_n\}$ is $C$-tight in $\D
([0,\tau])$ and then characterize the limit of every converging
subsequence. The $C$-tightness means that it satisfies the criteria for
tightness in the subspace $C$, as in
Theorem 11.6.3 of \cite{W02}.
The $C$-tightness implies that every subsequence has a further
converging subsequence with all limits having continuous sample paths w.p.1.
We demonstrate full convergence
by showing that all the convergent subsequences have the same limit.

\subsubsection{\texorpdfstring{Tightness of $\{W_n\}$}
{Tightness of \{W n\}}}\label{secWtightLLN}

First, the sequence $\{W_n\}$ is bounded, because $W_n (t) \ge0$ and
$W_n (t)$ increases at most at rate $1$.
The OL interval under question falls within a larger finite interval
$[0, \tau^*]$. Since the system has been assumed to start empty in
the initial UL interval,
a crude bound is $W_n (t) \le\tau^*$. Within the current OL interval,
we also can show that $W_n (0) \Rightarrow0$, so that $\limsup_{n \ra
\infty}{W_n (t)} \le\tau$.

Next, the modulus of continuity is bounded above because
$W_n (t + \delta) - W_n (t) \le\delta$.
It remains to bound $W_n (t + \delta) - W_n (t)$ below uniformly in $t$.
For that purpose, we work with the representation for $W_n$ in Section~\ref{secWrepPf}.
Let $\bar{I}_{n,j} \equiv n^{-1} I_{n,j}$ for $I_{n,j}$ in (\ref{Idecomp}).
We first observe that $n^{-1} I_{n,1} \Rightarrow0e$ and $n^{-1}
I_{n,2} \Rightarrow0e$ in $\D([0,\tau])$,
so that $\|\bar{I}_n - \bar{I}_{n,3}\|_{\tau} \Rightarrow0$.
However, by (\ref{Elim}), we already know that $\bar{E}_n \Rightarrow
E$ for $E(t) \equiv D(t) + s(t) - s(0)$.
Hence, we have $\bar{I}_{n,3} (t, \delta) - \bar{A}^I_n (t, \delta)
\Rightarrow E (t, \delta)$ in $\D([0,\tau])$
for $\bar{I}_{n,3}$ in (\ref{Idecomp}) and $\bar{A}^I_{n} (t, \delta
)$ in (\ref{in1}). However,\vspace*{1pt} by (\ref{Jbd3}), we can henceforth ignore
$\bar{A}^I_{n} (t, \delta)$.

By the assumptions for $\lambda$ and $F$ in Section~\ref{secFluid},
the integrand of $I_{n,3}$ in (\ref{Idecomp}) is bounded below by $c
\equiv F^{c} (\tau) \lambda_{\mathrm{inf}} (\tau) > 0$. Hence, we have the inequality
%
%
\begin{equation}
\label{limWrep} \frac{\bar{I}_{n,3}(t,\delta)}{c} \ge W_{n}(t) -W_{n}(t+
\delta) + \delta,
\end{equation}
so that we can write
%
%
\begin{equation}
\label{wBd1} W_{n}(t) -W_{n}(t+\delta) \le
\frac{\bar{I}_{n,3}(t,\delta)}{c} - \delta \le\frac{\bar{E}_n (t, \delta)}{c},
\end{equation}
and then combine the relations above to obtain
%
%
\begin{equation}
\label{wBd2} \limsup_{n \ra\infty}{\bigl\{W_{n}(t)
-W_{n}(t+\delta)\bigr\}} \le\frac
{D(t,\delta) + s(t, \delta)}{c} \equiv C \delta
\end{equation}
for some constant $C$.
Hence, the sequence $\{W_n (t)\dvtx  0 \le t \le\tau\}$ is $C$-tight. In
addition, the limit of any subsequence must be Lipschitz continuous.
Along the way, we have also shown that the sequences $\{\bar{I}_{n,3}
(t)\}$, $\{\bar{I}_{n} (t)\}$ and $\{\bar{A}^I_n (t)\}$ are tight as well.

\subsubsection{\texorpdfstring{Limit of convergent subsequences of $\{W_n\}$}
{Limit of convergent subsequences of \{W n\}}}
\label{secWsubLLN}

Since tightness implies that every subsequence has a convergent subsequence,
we complete the proof of the FWLLN for $W_n (t)$ by showing that every
convergent subsequence of $\{W_n\}$ converges to $w$ in $\D$.
It suffices to show that any limit of a convergent subsequence must satisfy
the ODE in (\ref{wODE}) w.p.1 or, equivalently, the integral
representation in (\ref{wIntEq}), because $w$
has been characterized as the unique solution to those equations.


First, by (\ref{enter1}) and (\ref{Elim}), we know that
%
%
\begin{equation}
\label{fwllnDview} \bar{E}_n (t,\varepsilon) \Rightarrow E(t, \varepsilon)
= \int_{t}^{t+
\varepsilon} b(s,0) \,ds
\end{equation}
in $\D$ as $n \ra\infty$. Moreover, as $\varepsilon\ra0$, the limit
in (\ref{fwllnDview}) approaches $b(t,0) = s(t)\mu+ \dot{s}(t)$.

We also consider the flow out of the queue in (\ref{in1}). Recall that
$\bar{I}_n (t, \ep)$ is asymptotically equivalent to $\bar{I}_{n,3}
(t, \ep)$ in (\ref{Idecomp}).
By the assumed convergence of $W_n \Rightarrow W$ and the continuous
mapping theorem applied to $\bar{I}_{n,3} (t, \ep)$, we have
%
%
\begin{equation}
\label{InLim}\qquad \bar{I}_n (t, \varepsilon) \Rightarrow I (t, \varepsilon)
\equiv\int_{t
- W(t)}^{t + \varepsilon- W( t+\varepsilon)} F^c (t-s)
\lambda(s) \,ds \qquad\mbox{in }\D\bigl([0, \tau]\bigr);
\end{equation}
that is, the limit $I (t, \ep)$ is determined once we know the limit
$W$. From (\ref{fwllnDview}) and (\ref{InLim}),
we also have
%
%
\begin{equation}
\label{AInLim} \quad\bar{A}^I_n (t, \varepsilon) =
\bar{E}_n (t, \varepsilon) - \bar{I}_n (t, \varepsilon)
\quad\Rightarrow\quad A^I (t, \ep) \equiv E (t, \ep) - I (t, \ep).
\end{equation}
Thus both limits $I (t, \ep)$ and $A^I (t, \ep)$ are determined given
the limit $W$.

In summary, we have the limits related by
%
%
\begin{eqnarray}
\label{fluidLimW} E(t, \varepsilon) & = & I(t, \varepsilon) + A^I (t,
\varepsilon)
\nonumber\\[-8pt]\\[-8pt]
& = & \int_{t - W(t)}^{t + \varepsilon- W(t+\varepsilon)} F^c(t-s)\lambda
(s) \,ds + A^I (t, \varepsilon).\nonumber
\end{eqnarray}
Again we can apply (\ref{Jbd3}) to deduce that $A^I (t, \ep)$ is
negligible relative to $I(t, \ep)$ for all suitably small $\ep$, so that
we can disregard $A^I (t, \ep)$ in (\ref{fluidLimW}).
Hence, combining (\ref{fwllnDview}), (\ref{fluidLimW}) and (\ref{Jbd3}),
we obtain
%
%
\begin{equation}
\label{Wconclusion}\qquad E(t, \varepsilon) = b(t, 0) \varepsilon+ o(\varepsilon) =
F^c\bigl(t-W(t)\bigr)\lambda \bigl(W(t)\bigr) \bigl(1 - \dot{W} (t)
\bigr) + o(\varepsilon)
\end{equation}
as $\varepsilon\downarrow0$ for almost all $t$ and almost all sample
paths of the limiting stochastic process $W$. In other words,
the proof of Theorem 3 of \cite{LW10} can be applied to $W$ to show
that $W$ satisfies the ODE (\ref{wODE}) w.p.1, that is, that Theorem 3 of
\cite{LW10} holds for $W$ w.p.1. Since there exists a unique solution
to that ODE, we must have $P(W = w) = 1$. Since this same conclusion
holds for all limits of
convergent subsequences,
we conclude that indeed $W_n \Rightarrow w$, as claimed.
Along the way, since we must have $W = w$, we determine the fluid
limits $I (t, \ep)$ and $A^I (t, \ep)$ as well; they are the limits
above with $W(t)$ replaced by $w(t)$.
We thus have two representations for $E (t) \equiv E(0,t)$,
%
%
\begin{equation}
\label{ElimRep}\quad E(t) = D(t) + s(t) - s(0) = \int_{0}^{t-w(t)}
F^c(t-s)\lambda(s) \bigl(1 - \dot{w} (s)\bigr) \,ds.
\end{equation}

\subsubsection{\texorpdfstring{The FWLLN for $V_n(t)$}
{The FWLLN for V n(t)}}\label{secVfwlln}

By the definitions of the HWT and PWT, we necessarily have
the PWT $V_n$ satisfying the equation
%
%
\begin{eqnarray}
\label{Vrep} V_n \bigl(t - W_n (t)\bigr) & = &
W_n (t) + O(1/n) \quad\mbox{or, equivalently}
\nonumber\\[-8pt]\\[-8pt]
V_n (t) & = & W_n \bigl(t + V_n (t)+
O(1/n)\bigr)+ O(1/n),\nonumber
\end{eqnarray}
given $W_n (t)$. Note that these equations relating the PWT and HWT for
the stochastic queueing systems are slightly different from those for
the deterministic fluid models as in (\ref{vEq}). The first equation
in (\ref{Vrep}) holds since the PWT at $t-W_{n}(t)$ equals the HWT at
$t$ plus the remaining time until the first busy server becomes
available, which is $O(1/n)$ since there are $O(n)$ servers. The second
equation in (\ref{Vrep}) holds simply by applying a change of variable
in the first equation.

We already have established the FWLLN for $W_n (t)$, yielding $W_n
\Rightarrow w$, where $w$ is a continuous function.
Moreover, $w$ has left and right derivatives everywhere, which are bounded.

We now exploit Theorems 3--6 of \cite{LW10} establishing key
properties of the HWT and PWT fluid functions $w$ and $v$.
The additional property (\ref{flowInto}) here implies that there exists
a constant $\gamma> 0$ such that $\dot{w}(t) < 1 - \gamma$, $0 \le t
\le\tau$.
By Theorems 5 and 6 of \cite{LW10}, $v$ is continuous, where $v$ is
the unique solution to the corresponding fluid equations, for example,
as in (\ref{vEq}).
Hence, from the construction of $v$ in the proof of Theorem 5 of \cite
{LW10} and (\ref{vDeriv}),
we deduce, first, for the given fluid functions $(w, v)$ and any other
$(w_1, v_1)$ that
$\| v_1 - v\|_{\tau} < \| w_1 - w\|_{\tau}/\gamma$.
Hence, we deduce that
%
%
\begin{equation}
\label{normBd} \| V_n - v\|_{\tau} < \| W_n -
w\|_{\tau}/\gamma O(1/n).
\end{equation}
Since, $\| W_n - w\|_{\tau} \Rightarrow0$, also $\| V_n - v\|_{\tau}
\Rightarrow0$.
Hence, the proof of the FWLLN is complete.

\subsection{\texorpdfstring{Proof of Theorem \protect\ref{thFCLT}: The FCLT}
{Proof of Theorem 4.2: The FCLT}}\label{fcltPf}

We now turn to the proof of the FCLT, still under our special initial
conditions imposed in Section~\ref{secSpecialInitial}.
From Section~\ref{secQlimPf}, we know that, for the queue length $Q_n
(t)$ and the number in system $X_n (t)$,
it suffices to prove convergence of the scaled waiting times $\hat{W}_n$.
Just as for the FWLLN, we do this in two steps. We first prove
tightness and then we characterize the limit of all convergent subsequences.

\subsubsection{\texorpdfstring{Tightness of the sequence $\{\hat{W}_n\}$}
{Tightness of the sequence \{W n\}}}\label{secWhatTightPf}

We start by proving $C$-tightness of the sequence
$\{\hat{W}_n\} \equiv\{\sqrt{n} (W_n (t) - w(t))\}$.
We do a proof
by contradiction.
First, suppose that $\{\hat{W}_n\}$ is not stochastically bounded;
that is, for all real numbers $M > 0$ no matter how large and for all
$\varepsilon>0$
no matter how small,
there exists
$n$ such that $P(\|\hat{W}_n\|_{\tau} > M) > \varepsilon$.
However, from Sections~\ref{secWrepPf} and \ref{secWsubLLN},
including (\ref{ElimRep}), we know that $\|\hat{E}_n - \hat{I}_n\|
_{\tau} \Rightarrow0$,
where
%
%
\begin{equation}
\label{EWbd} \hat{I}_n (t) \equiv\sqrt{n}\int_{t - w(t)}^{t-W_n (t)}
F^c(t-s)\lambda(s) \bigl(1 - \dot{w} (s)\bigr) \,ds.
\end{equation}
Hence, there exists $n$ for all $M > 0$, no matter how large and for
all $\varepsilon>0$
no matter how small, such that
%
%
\begin{equation}
\label{EWbd2} P\bigl(\|\hat{E}_n \|_{\tau} \ge c \|
\hat{W}_n \|_{\tau} \ge cM\bigr) > \ep,
\end{equation}
where
$c$ is the strictly positive infimum of the integrand in (\ref{EWbd})
(because $\lambda(t) > \lambda_{\mathrm{inf}} > 0$, $w(t) < 1$
and $w$ is uniformly continuous on the interval $[0, \tau]$).
However, this would contradict the established convergence $\hat{E}_n
\Rightarrow\hat{E}$ in (\ref{Elim}) and (\ref{Elimit}).
Hence the sequence $\{\hat{W}_n\}$ must actually be stochastically bounded.

Second, even though the sequence $\{\hat{W}_n\}$ is
stochastically bounded, it is possible that the modulus of $\{\hat
{W}_n\}$ is
not asymptotically negligible, as in (11.6.4) of \cite{W02}.
Thus, suppose that there exists $\varepsilon> 0$ and $\delta> 0$
such that
%
%
\begin{equation}
\label{modulus} P \bigl(\omega_{\hat{W}_n} (\delta) > \varepsilon\bigr) > \delta
\end{equation}
for all $\delta> 0$, no matter how small,
and some $n$, no matter how large, where
%
%
\begin{equation}
\label{modulus2} \omega_{x} (\delta) \equiv\sup_{0 \le t < t_1 < t_2 \le t + \delta
\le\tau}{
\bigl\{\bigl|x(t_2) -x(t_1)\bigr|\bigr\}}.
\end{equation}
Consider a subsequence of
$n$ for which this is true. Then there must exist a sequence
$\{(\delta_n, t_n)\}$ where $0 \le t_n < \tau$ and $\delta_n
\downarrow0$ as $n \ra\infty$ such that $P(|\hat{W}_n (t_n +
\delta_n) - \hat{W}_n (t_n)| > \gamma) > \varepsilon$ for all $n$.
Since, $0 \le t_n \le\tau$ for all $n$, there exists a
convergent subsequence of $\{t_n\}$. So it suffices to assume
that $t_n \ra t$ as $n \ra\infty$.

We now work with $I_{n,3}$ in (\ref{Idecomp}), using the fluid limits
$I$ and $A^I$ determined in Section~\ref{secWsubLLN}, that is,
%
%
\begin{equation}
\label{fkuidIA}\qquad I(t, \delta) = \int_{t - w(t)}^{t + \delta- w( t+\delta)}
F^c(t-s)\lambda(s) \,ds \quad\mbox{and}\quad A^I (t) = E(t) - I(t).
\end{equation}
Thus, by the continuity of $\tilde{q}$, the $\sqrt{n}$-scaled process
satisfies
%
%
\begin{eqnarray}
&& \hat{I}_{n,3} (t_n,\delta_n)\nonumber\hspace*{-35pt}\\
&&\quad = \sqrt{n} \biggl(\int_{t_n - W_n (t_n)}^{t_n + \delta_n - W_n
(t_n + \delta_n)} \tilde{q}(t,t-s) \,ds - \int_{t_n - w (t_n)}^{t_n + \delta_n -
w (t_n + \delta_n)} \tilde {q}(t, t-s) \,ds \biggr)\nonumber\hspace*{-35pt}\\
&&\quad = \sqrt{n} \biggl(\int_{t_n + \delta_n - w (t_n + \delta_n)}^{t_n +
\delta_n - W_n (t_n + \delta_n)} \tilde{q}(t,t-s) \,ds - \int_{t_n - W (t_n)}^{t_n -
w (t_n)} \tilde{q}(t, t-s) \,ds \biggr)\hspace*{-35pt}\\
&&\quad =\tilde{q}\bigl(t_n + \delta_n, w(t_n +\delta_n)\bigr) \hat{W}_n (t_n +
\delta_n) - \tilde{q}\bigl(t_n, w(t_n)\bigr) \hat{W}_n (t_n) + o(1)\nonumber\hspace*{-35pt}\\
&&\quad = \tilde{q}\bigl(t, w(t)\bigr) \bigl(\hat{W}_n(t_n + \delta_n) - \hat {W}_n
(t_n)\bigr) + o(1) \qquad\mbox{as } n \ra\infty,\nonumber\hspace*{-35pt}
\end{eqnarray}
so that
%
%
\begin{equation}
\label{hatIn222}
\quad \limsup_{n \ra\infty}{\bigl|\hat{I}_{n,3}
(t_n, \delta_n)\bigr|} \ge\tilde {q}\bigl(t, w(t)\bigr)
\limsup_{n \ra\infty}{\bigl| \hat{W}_n (t_n +
\delta_n) - \hat{W}_n (t_n)\bigr|}.
\end{equation}
Since limits have been established for the sequences $\hat{I}_{n,1}$
and $\hat{I}_{n,2}$, (\ref{hatIn222}) implies that, for some $\gamma'
> 0$,
%
%
\begin{eqnarray}
\label{hatIn3Tight} && \limsup_{n \ra\infty}{P\bigl(\bigl|
\hat{I}_{n} (t_n, \delta_n)\bigr| > \gamma
'\bigr)}
\nonumber\\[-8pt]\\[-8pt]
&&\qquad \ge \limsup_{n \ra\infty}{P\bigl(\bigl| \hat{W}_n
(t_n + \delta_n) - \hat{W}_n
(t_n)\bigr| > \gamma\bigr)} > 0.\nonumber
\end{eqnarray}

However, together with (\ref{Jbd3}), which implies that $|\hat{A}^I_n
(t, \ep)| \le K \ep|\hat{I}_n (t, \ep)|$
for some constant $K$ for all $\ep$ suitably small, uniformly in $n$
and $t$, the limit in (\ref{hatIn3}) implies
that we cannot have $\hat{E}_n \Rightarrow\hat{E}$ as indicated in
(\ref{Elim}),
which is a contradiction.
Hence, the modulus property for the sequence $\{\hat{W}_n\}$ in (\ref
{modulus}) must actually not hold.
Thus, we have shown that the sequence $\{\hat{W}_n\}$ must in fact be tight.

\subsubsection{\texorpdfstring{Characterizing the limit $\hat{W}$}
{Characterizing the limit W}}\label{secWcharPf}

We now characterize the limit of any convergent subsequence of
the sequence $\{\hat{W}_n\}$. Without changing the notation, suppose
that $\hat{W}_n \Rightarrow\hat{W}$
through some subsequence. Of course, we also have $W_n \Rightarrow w$
along this subsequence and all the other fluid limits.
We thus want to characterize the distribution of $\hat{W}$.
To do so, we again exploit the representation of the flow into service, writing
%
%
\begin{equation}
\label{r1} E_n (t) = \sum_{i = 1}^{\lfloor t/\ep\rfloor}
E_n\bigl((i-1)\varepsilon, \varepsilon\bigr) + E_{n,r} (t,
\varepsilon),
\end{equation}
where $E_{n,r} (t, \varepsilon)$ is the final remainder term associated
with a final partial interval
and
$E_n (t, \varepsilon) = I_n (t, \varepsilon) - A^I_n (t, \varepsilon)$ as in
(\ref{in1})
with $I_n (t, \varepsilon) = I_{n,1} (t, \varepsilon) + I_{n,2} (t,
\varepsilon) + I_{n,3} (t, \varepsilon)$ as in (\ref{in3})
and (\ref{Idecomp}).
Since we have established that\vadjust{\goodbreak} $\hat{E}_n (t) \Rightarrow\hat{E}
(t)$ in $\D([0,\tau])$, as stated in (\ref{Elim}) and (\ref{Elimit}),
we can ignore the final remainder term in (\ref{r1}). The $C$-tightness
following from the convergence implies that the scaled remainder
term is asymptotically negligible.

For any $t > 0$ (which applies to $i \ge1$), let the $\sqrt {n}$-scaled processes over the intervals $[t, t + \varepsilon]$ be
%
%
\begin{equation}
\label{r2} \hat{I}_{n,j} (t, \varepsilon) \equiv n^{-1/2}
\bigl(I_{n,j} (t, \varepsilon) - n I_{j} (t, \varepsilon)\bigr),
\end{equation}
where $I_j (t, \varepsilon)$ has been determined, and similarly for the
other processes.
In Section~\ref{secWtightLLN} we observed that $I_1 = I_2 = 0e$.

By (\ref{Idecomp}), the FWLLN for $W_n$ and the FCLT for
$\hat{Z}_{n,1}$ in (\ref{tildeLim}),
%
%
\begin{eqnarray}
\label{hatIn1} &&\hat{I}_{n,1} (t, \ep) \equiv \frac{1}{\sqrt {n}}I_{n,1}(t,
\varepsilon) \nonumber\\[-8pt]\\[-8pt]
&&\quad\Rightarrow\quad \hat{I}_1 (t, \ep)
\equiv  c_{\lambda} \int_{t-w(t)}^{t+\varepsilon-w(t+\varepsilon
)}F^c(t-s)\,d
\tilde{\sB}_{\lambda}\bigl(\Lambda(s)\bigr),\nonumber
\end{eqnarray}
where
$\sB_{\lambda}$ is the BM associated with the arrival process, and
$c_{\lambda}^2$ is its variability parameter, as in (\ref{arriveFCLT}).

Similarly, by (\ref{Idecomp}), the FWLLN for $W_n$ and the FCLT for
$\hat{Z}_{n,2}$ in (\ref{tildeLim}),
%
%
\begin{eqnarray}
\label{hatIn2} &&\hat{I}_{n,2} (t, \varepsilon) \equiv
\frac{1}{\sqrt {n}}I_{n,2}(t,\varepsilon) \nonumber\\
&&\quad\Rightarrow\quad \hat{I}_{2}
(t, \varepsilon)
\equiv \int_{t-w(t)}^{t+\varepsilon-w(t+\varepsilon)}\int_{0}^{\infty}{
\mathbf1}(s+x>t)\,d \sR(s,x)
\\
&&\hphantom{\quad\Rightarrow\quad \hat{I}_{2}
(t, \varepsilon)}
\deq -\int_{t-w(t)}^{t+\varepsilon-w(t+\varepsilon)}\sqrt {F(t-s)F^c
(t-s)}\,d \tilde{\sB}_{a}\bigl(\Lambda(s)\bigr),\nonumber
\end{eqnarray}
where
$\sB_{a}(\cdot)$ is a BM associated with the patience times.

For $I_{n,3}$ in (\ref{Idecomp}), we first write
%
%
\begin{equation}
\label{E3hat} \hat{I}_{n, 3} (t,\varepsilon) \equiv n^{-1/2} \bigl(
I_{n,3} (t,\varepsilon) - n I_3 (t, \ep) \bigr)
\end{equation}
for $I_3$ previously determined, that is,
%
%
\begin{equation}
\label{fluidE3} I_3 (t, \ep) \equiv\int_{t-w(t)}^{t+\varepsilon-w(t+\varepsilon
)}F^{c}(t-s)
\lambda(s)\,ds.
\end{equation}
Then, exploiting the assumed convergence $\hat{W}_n \Rightarrow\hat
{W}$ along the subsequence, we obtain
%
%
\begin{eqnarray}
\label{hatIn3}
\quad&&\hat{I}_{n,3}(t, \varepsilon) \nonumber\\
&&\qquad = \frac{1}{\sqrt{n}}
\biggl(n\int_{t-W_{n}(t)}^{t+\varepsilon-W_{n}(t+\varepsilon)}\tilde{q}(t,t-s)\,ds - n\int
_{t-w(t)}^{t+\varepsilon-w(t+\varepsilon)}\tilde{q}(t,t-s)\,ds \biggr)
\nonumber
\\
&&\qquad = \sqrt{n} \biggl(\int_{t-W_{n}(t)}^{t-w(t)}
\tilde{q}(t,t-s)\,ds + \int_{t+\varepsilon-w(t+\varepsilon)}^{t+\varepsilon-W_{n}(t+\varepsilon
)}
\tilde{q}(t,t-s)\,ds \biggr)
\nonumber\\[-8pt]\\[-8pt]
&&\qquad = \tilde{q}\bigl(t,w(t)\bigr) \sqrt{n}\bigl(W_{n}(t)-w(t)\bigr)
\nonumber\\
&&\qquad\quad{} - q\bigl(t,w(t+\varepsilon)-\varepsilon\bigr) \sqrt{n}\bigl(W_{n}(t+
\varepsilon )-w(t+\varepsilon)\bigr) + o(1)
\nonumber
\\
&&\quad \Rightarrow\quad \hat{I}_3 (t, \ep) \equiv\tilde{q}\bigl(t,w(t)\bigr)
\hat {W} (t) - \tilde{q}\bigl(t,w(t+\varepsilon)-\varepsilon\bigr) \hat{W} (t +
\varepsilon)\nonumber
\end{eqnarray}
as $n \ra\infty$. Exploiting (\ref{Jbd3}), we see that $\hat{A}^I_n$
is asymptotically negligible compared to $\hat{I}_{n,s}$. Hence,
we have established the convergence
%
%
\begin{equation}
\label{vectorLim} \bigl(\hat{W}_n, \hat{E}_n (t, \ep),
\hat{I}_n (t, \ep)\bigr) \Rightarrow \bigl(\hat{W}, \hat{E} (t, \ep),
\hat{I} (t, \ep)\bigr)
\end{equation}
in $\D([0, \tau]) \times\D^3([t, t+\ep])$, where
%
%
\begin{equation}
\label{KeyRel} \hat{E} (t, \ep) = \hat{I}_1 (t, \ep) +
\hat{I}_2 (t, \ep) + \hat {I}_3 (t, \ep) \bigl(1 + o(\ep)
\bigr)
\end{equation}
with all the limits having been identified explicitly.
Substituting the established limits into (\ref{KeyRel}), we obtain
%
%
\begin{eqnarray}
\label{preSDE}
&&\sB_s \bigl(D(t+ \ep)\bigr) - \sB_s
\bigl(D(t)\bigr) \nonumber\\
&&\qquad= \int_{t-w(t)}^{t+\varepsilon
-w(t+\varepsilon)}F^c(t-s)\,d
\bigl(c_{\lambda} \tilde{\sB}_{\lambda}\bigl(\Lambda (s)\bigr)\bigr)
\nonumber\\[-8pt]\\[-8pt]
&&\qquad\quad{} -\int_{t-w(t)}^{t+\varepsilon-w(t+\varepsilon)}\sqrt {F(t-s)F^c(t-s)}\,d
\tilde{\sB}_{a}\bigl(\Lambda(s)\bigr)
\nonumber\\
&&\qquad\quad{} + \tilde{q}\bigl(t,w(t)\bigr) \hat{W} (t) - \tilde{q}\bigl(t,w(t+\varepsilon )-
\varepsilon\bigr) \hat{W} (t + \varepsilon) + o(\ep) \qquad\mbox{as }\ep\downarrow0.
\nonumber
\end{eqnarray}
Moreover, for each $\ep> 0$, we have the corresponding limit for the
sum $\hat{E}_n (t)$ in (\ref{r1}).
As $\ep\downarrow0$, this sum converges
in mean square to the stochastic integral associated with a stochastic
differential equation (SDE)
determined by (\ref{preSDE}). Thus, the distribution of $\hat{W}$ is
determined by this SDE.
The SDE is well defined because all but the term $\hat{W} (t + \ep) -
\hat{W(t)}$ involve
BM terms, which produce known differential terms.
In particular, using informal differential notation, we see that, as
$\ep\downarrow0$,
\begin{eqnarray*}
\hat{E} (t, \ep) & \ra& d \tilde{\sB}_s \bigl(D(t)\bigr),
\\
\hat{I}_1 (t, \ep) & \ra& d\hat{I}_1 (t) \equiv
c_{\lambda} F^c \bigl(w(t)\bigr) \,d \tilde{\sB}_{\lambda} (
\Lambda\bigl(t - w(t)\bigr),
\\
\hat{I}_2 (t, \ep) & \ra& d\hat{I}_2 (t) \equiv- \sqrt{F
\bigl(w(t)\bigr) F^c \bigl(w(t)\bigr)} \,d \tilde{\sB}_a (
\Lambda\bigl(t - w(t)\bigr),
\\
\hat{I}_3 (t, \ep) & \ra& d\hat{I}_3 (t) \equiv-\tilde
{q}\bigl(t,w(t)\bigr)\,d\hat{W}(t) + \bigl(1-\dot{w}(t)\bigr)
\tilde{q}_{x}\bigl(t,w(t)\bigr)\,dt \hat{W}(t),
\end{eqnarray*}
where we exploit the assumed differentiability of the arrival rate
function $\lambda$ and
\begin{eqnarray*}
&& \frac{\tilde{q}(t,w(t))-\tilde{q}(t,w(t+\varepsilon)-\varepsilon
)}{\varepsilon}
\\
&&\qquad = \biggl(\frac{\tilde{q}(t,w(t))
-\tilde{q}(t,\varepsilon+w(t+\varepsilon))}{w(t)+\varepsilon-w(t+\varepsilon
)} \biggr) \biggl(\frac{w(t)+\varepsilon-w(t+\varepsilon)}{\varepsilon} \biggr)
\\
&&\qquad \rightarrow\tilde{q}_{x}\bigl(t,w(t)\bigr)\bigl[1-\dot{w}(t)
\bigr] \qquad\mbox{as } \varepsilon\rightarrow0
\end{eqnarray*}
in the treatment of $\hat{I}_3 (t, \ep)$.

Putting the $d\hat{W}(t)$ term on the left, and thus expressing it in
terms of all others, we get
the SDE
%
%
\begin{eqnarray}
\label{sdePf} d\hat{W}(t) & = & h(t) \hat{W}(t)\,dt - \biggl(\frac{1}{q(t,w(t))}
\biggr)\,d \tilde{\sB}_{s}\bigl(D(t)\bigr)
\nonumber
\\
&&{} - \biggl(\frac{\sqrt{F(w(t))F^c (w(t))}}{q(t,w(t))} \biggr)\,d \tilde{\sB}_{a}\bigl(
\Lambda\bigl(t-w(t)\bigr)\bigr)
\nonumber
\\
&&{} + \frac{F^c(w(t))c_{\lambda}}{q(t,w(t))}\,d \tilde{\sB }_{\lambda}\bigl(\Lambda\bigl(t-w(t)
\bigr)\bigr)
\\
& = & h(t) \hat{W}(t)\,dt +I_{s}(t) \,d \sB_{s}(t) +
I_{a}(t)\,d \sB _{a}(t) + I_{\lambda}(t)\,d
\sB_{\lambda}(t)
\nonumber
\\
& = & h(t)\hat{W}(t)\,dt +I(t)\,d \sB(t)\nonumber
\end{eqnarray}
as in (\ref{sde}), where $h(t)$, $I_1 \equiv I_{\lambda}$, $I_2
\equiv
I_{s}$ and $I_s\equiv I_{a}$ are given in (\ref{kernels}) and $I$ is
given in (\ref{Isquared}).
while
$\sB_{1} \equiv\sB_{\lambda}$, $\sB_{2} \equiv\sB_{s}$ and $\sB
_{3} \equiv\sB_{a}$ are all independent standard BMs.

We claim that the SDE in (\ref{sdePf}) and (\ref{sde}) has the
analytic solution
%
%
\begin{eqnarray}
\hat{W}(t) & = & \hat{W}(0) H(t,0) + \int_{0}^{t}
\biggl( - \frac
{1}{q(u,w(u))} \biggr) H(t,u)\,d \tilde{\sB}_{s}
\bigl(D(u)\bigr)
\nonumber
\\
&&{} + \int_{0}^{t} \biggl( - \frac{\sqrt{F(w(u))
F^c(w(u))}}{q(u,w(u))}
\biggr) H(t,u)\,d \tilde{\sB}_{a}\bigl(\Lambda \bigl(u-w(u)\bigr)\bigr)
\nonumber
\\
\label{WB}
&&{} + \int_{0}^{t} \frac{F^c(w(u))c_{\lambda}}{q(u,w(u))} H(t,u) \,d
\tilde{\sB}_{\lambda}\bigl(\Lambda\bigl(u-w(u)\bigr)\bigr)
\\
\label{SDEsol}
& \deq& \hat{W}(0) H(t,0) + \int_{0}^{t}H(t,u)
I(u)\,d \sB(u)
\\
& \deq& \hat{W}_{0}(t) +\hat{W}_{\lambda}(t) +
\hat{W}_{s}(t) + \hat{W}_{a}(t),
\nonumber
\end{eqnarray}
where $\hat{W}_{0} = 0e$, $\hat{W}_{1} \equiv\hat{W}_{\lambda}$,
$\hat{W}_{2} \equiv\hat{W}_{s}$ and $\hat{W}_{3} \equiv\hat
{W}_{a}$ are independent processes,
as given in Theorem \ref{thFCLT}.

We verify (\ref{SDEsol}) from (\ref{sde}) using It\^{o}'s formula. Let
$f(x,t)\equiv x e^{-\int_{0}^{t}h(v)\,dv}$, we have by It\^{o}'s formula that
\begin{eqnarray*}
df\bigl(\hat{W}(t),t\bigr) & = & e^{-\int_{0}^{t}h(v)\,dv} \,d\hat{W}(t) - h(t)
e^{-\int_{0}^{t}h(v)\,dv} \hat{W}(t)\,dt,
\\
& = & e^{-\int_{0}^{t}h(v)\,dv} I(t)\,d \sB(t).
\end{eqnarray*}
Integrating both sides yields
\[
e^{-\int_{0}^{t}h(v)\,dv} \hat{W}(t) = \hat{W}(0) + \int_{0}^{t}e^{-\int_{0}^{u}h(v)\,dv}
I(u)\,d \sB(u),
\]
from which (\ref{SDEsol}) follows by multiplying through by $H(t,0)
\equiv e^{\int_{0}^{t}h(v)\,dv}$.

\subsubsection{FCLT for other processes}\label{secFCLTother}

So far, we have established the FCLT for the HWT process $W_n (t)$,
still under the special initial condition starting with all servers
busy and an empty queue.
We now use this result to establish limits for the other processes,
under this same initial condition.\vspace*{9pt}

\noindent\textit{The queue length and the number in system.\quad}
We now obtain the limit for $\hat{Q}_n$ and $\hat{X}_n$ from (\ref
{Qlim}) and (\ref{XQdiff}), incorporating the limit for $\hat{W}_n$ into
$\hat{Q}_{n,3}$,
using the limit $\hat{W}_n \Rightarrow\hat{W}$ just established. We
obtain the expression in Theorem \ref{thFCLT}
by putting the contributions from the arrival process, service times
and patience times into their respective terms.
We have thus established the FWLLN in Theorem \ref{thFWLLN} and the
FCLT in Theorem \ref{thFCLT} under the special initial condition,
in which all servers are busy, and the queue is empty at time $0$, the
beginning of the OL interval.\vspace*{9pt}

\noindent\textit{The potential waiting time.\quad}
We start with the fluid equation $v(t) = w(t + v(t))$ in (\ref{vEq})
and the corresponding equation for the queueing models,
$V_n(t) = W_n (t + V_n (t)+O(1/n))+O(1/n)$, as in (\ref{Vrep}). Let
$\Delta V_n (t) \equiv V_n (t) - v(t)$
and $\Delta W_n (t) \equiv W_n (t) - w(t)$.
We exploit the differentiability of $w(t)$ with $\dot{w} (t) < 1 - \ep
$ for some $\ep> 0$,
the differentiability of $\dot{w}$ [because we assumed that $\lambda$
is differentiable in order to have $\tilde{q}_x (t,x)$ well defined] and
Taylor's theorem to write
%
%
\begin{eqnarray}
\label{v101}\quad \Delta V_n (t) & = & \Delta W_n \bigl(t +
V_n (t)+O(1/n)\bigr) + w\bigl(t + V_n (t)\bigr)- w
\bigl(t + v(t)\bigr)+O(1/n)
\nonumber\\
& = & \Delta W_n \bigl(t + V_n (t)+O(1/n)\bigr) +
\dot{w}\bigl(t + v(t)\bigr) \Delta V_n (t)
\\
&&{} + \ddot{w}\bigl(t + v(t)\bigr) \frac{(\Delta V_n (t))^2}{2} + o\bigl(\bigl(\Delta
V_n (t)\bigr)^2\bigr)+O(1/n).\nonumber
\end{eqnarray}
We exploit the FCLT for $W_n (t)$, the FWLLN for $V_n (t)$ and the
continuous mapping theorem
to get $\sqrt{n}\Delta W_n (t + V_n (t)) = \hat{W}_n (t + V_n (t))
\Rightarrow\hat{W}(t + v(t))$ in $\D([0, \tau])$.
From (\ref{normBd}), we see that there exists $\gamma> 0$ such that
%
%
\begin{equation}
\label{v102} \| \Delta V_n\|_{\tau} \le\frac{ \|\Delta W_n \|_{\tau}}{\gamma} +
O(1/\sqrt{n})= O(1/\sqrt{n}) \qquad\mbox{as } n \ra\infty.
\end{equation}
We can then apply (\ref{v102}) with the two-term expansion in (\ref
{v101}) to get
%
%
\begin{equation}
\label{v103} \sup_{0 \le t \le\tau}{ \biggl\{\biggl\llvert
\hat{V}_n (t) - \frac{\hat
{W}_n (t + v(t))}{1 - \dot{w} (t + v(t))} \biggr\rrvert \biggr\} } = \sqrt
{n} O\bigl(\bigl(\|\Delta V_n \|_{\tau}\bigr)^2\bigr) =
O(1/\sqrt{n}),\hspace*{-35pt}
\end{equation}
proving that
%
%
\begin{equation}
\label{v104} \hat{V}_n (t) \Rightarrow\hat{V} (t) \equiv
\frac{\hat{W} (t +
v(t))}{1 - \dot{w} (t + v(t))} \qquad\mbox{in }\D\bigl([0,\tau]\bigr)
\end{equation}
as claimed.\vspace*{9pt}

\noindent\textit{The abandonment process $A_n (t)$.\quad}
We obtain the limits for $\bar{A}_n$ and $\hat{A}_n$ in (\ref{fLim})
and (\ref{gXlim}) directly from the flow conservation
representation in (\ref{aban})
and the established limits above. We see that $\bar{A}_n \Rightarrow
A$ in $\D$
and $\hat{A}_n \Rightarrow\hat{A}$ in $\D$, jointly with the other
processes, for $\hat{A}$ in (\ref{Glims}).

\subsection{\texorpdfstring{Treating the initial conditions in (\protect\ref{InitialLim})}
{Treating the initial conditions in (4.1)}}\label{secInitialPf}

It now remains to extend the FWLLN and the FCLT for the number in
system in an OL interval to the general initial condition given in
(\ref{InitialLim}).
As in the statement of Theorem \ref{thFCLT}, let $X_n (t)$ be the
number in system during the OL interval with the initial condition
(\ref{InitialLim}), and let
$X^{*}_n (t)$ be the number in system during the OL interval starting
with all servers busy and an empty queue,
for which we have proved the FWLLN and FCLT in the preceding subsections.

We assume that the two processes $X_n (t)$ and $X^{*}_n (t)$ are
defined on the same probability space,
having the same arrival process, service times and abandonment times,
with the service times and abandonment times assigned in order of
customers entering service and the queue, respectively.
These processes differ by the initial conditions $X_n (0) - s_n (0)$,
for which the scaled versions
have been assumed to converge in (\ref{InitialLim}). However, we need to
carefully consider the consequence of this difference at time $0$ as
time evolves
within the interval $[0, \tau]$.

We establish the desired limits for $X_n (t)$ by showing that
%
%
\begin{equation}
\label{asymEq2} \bigl\|\bar{X}_n - \bar{X}^{*}_n
\bigr\|_{\tau} \Rightarrow0 \quad\mbox{and}\quad\bigl\|\hat {X}_n - \bigl(
\hat{X}^{*}_n + \hat{X}_n(0)
F^c_w (t)\bigr)\bigr\|_{\tau} \Rightarrow0
\end{equation}
in $D([0,\tau])$ as $n \ra\infty$, where $\hat{X}_n (0)$ is
independent of $\hat{X}^{*}_n$ and $F^c_w (t)$ is given in (\ref{Glims})
with $w(t)$ being the HWT in the fluid model and $h_F (x) \equiv
f(x)/F^c (x)$ being the hazard rate function of the patience c.d.f. $F$,
which is positive and bounded
by previous assumptions on $F$. As a consequence of the first limit in
(\ref{asymEq2}), the fluid limit appearing in the centering terms of
the scaled processes $\hat{X}^{*}_n$ and $\hat{X}_n$ are identical.

We now proceed to justify (\ref{asymEq2}).
Since the customers enter service in a FCFS order, the excess customers
at time $0$ soon enter service.
However, the excess still remains
because new customers\vadjust{\goodbreak} arrive and join the queue to replace those that
entered service.
An important insight is the observation that the remaining excess can
always be considered among those customers
that have been in the system for the longest time among all waiting customers.

Since the abandonment hazard rate
is bounded above, the abandonment rate is controlled.
Since the fluid model is in an OL interval with $\Lambda(t) > D(t)$
for all $t > 0$,
with the initial net input rate to service $\lambda(0) - s(0) \mu-
\dot{s} (0) > 0$,
the servers become all busy and remain so afterwards in an interval
$[t_{1,n}, t_2]$ for $0 < t_{1,n} = O(1/\sqrt{n}) < t_2 < \tau$.
Thus there are at most $O(\sqrt{n})$ empty servers for a period of
only $O(1/\sqrt{n})$.
Hence, the difference between $X_n^{*} (t)$ and $X_n (t)$ is
asymptotically only the initial difference
adjusted by abandonments over the interval $[0,t]$. In particular, we have
%
%
\begin{equation}
\label{Xdiff} \bigl\|X_n - \bigl(X^{*}_n +
U^{+}_n - U^{-}_n\bigr)
\bigr\|_{\tau} = O(1) \qquad\mbox{as } n \ra \infty,
\end{equation}
where
%
%
\begin{eqnarray}
\label{Udef} U^{+}_n (t) & \equiv& \bigl(X_n
(0) - s_n (0)\bigr)^+ - A_{i,n,+} (t),
\nonumber\\[-8pt]\\[-8pt]
U^{-}_n (t) & \equiv& -\bigl(X_n (0) -
s_n (0)\bigr)^- - A_{i,n,-} (t)\nonumber
\end{eqnarray}
with $(x)^- \equiv\min{\{x,0\}}$,
$A_{i,n,+} (t)$ being the number of abandonments from the initial
positive excess number of customers, $(X_n (0) - s_n (0))^+ > 0$,
given that it is positive, while $A_{i,n,-} (t)$ is the number of
abandonments from the positive difference $-(X_n (0) - s_n (0))^- $,
resulting from a initial negative excess number of customers, $(X_n (0)
- s_n (0))^- < 0$,
given that it is indeed negative.
Fortunately, the limiting behavior of $A_{i,n,+} (t)$ and $A_{i,n,-}
(t)$ are essentially the same,
so that we need not treat the positive part and the negative part differently.

We are now ready to prove the FWLLN.
Since
$0 \le(X_n (0) - s_n (0))^+ - A_{i,n,+} (t) \le(X_n (0) - s_n (0))^+
= O(\sqrt{n})$
and
$0 \le-(X_n (0) - s_n (0))^-  -\break A_{i,n,-} (t) \le-(X_n (0) - s_n
(0))^- = O(\sqrt{n})$,
we deduce that $\|\bar{X}_n - \bar{X}^{*}_n\|_{\tau} \Rightarrow0$
as $n \ra\infty$.
Hence, we have completed proof of the FWLLN $\bar{X}_n \Rightarrow X$
in $\D([0, \tau])$.
The rest of Theorem \ref{thFWLLN} follows for the general initial
conditions $\bar{X}_n (0) \Rightarrow X (0)$ as well.

We now turn to the FCLT. We will show that $\|\hat{X}_n - (\hat
{X}^{*}_n + \hat{X}_n (0) F^c_w (\cdot))\|_{\tau} \Rightarrow0$
in $D([0,\tau])$, as in (\ref{asymEq2}).
For that, we need to carefully examine the processes $A_{i,n,+} (t)$
and $A_{i,n,-} (t)$,
recording the number of abandonments from the deviation $X_n (0) - s_n (0)$.
Suppose that $X_n (0) - s_n (0) > 0$,
so that we focus on $A_{i,n,+} (t)$.
Since the abandonments $A_{i,n,+} (t)$ always come from the waiting
customers that have been in the system the longest,
which means at the right boundary of the queue length process, which
asymptotically is at $w(t)$,
the abandonment making up $A_{i,n,+} (t)$ occurs asymptotically at rate
$h_F (w (u))$ at time $u$ through all time.

Of course, specific abandonments are random. Nevertheless, because the
size of the deficiency is order $O(\sqrt{n})$
and we scale by dividing by $\sqrt{n}$ when we scale for the FCLT, the
impact actually becomes deterministic, by the FWLLN (or
Glivenko--Cantelli theorem).
In particular, the FCLT-scaled version of the process $U_n^+ (t)$ in
(\ref{Udef}) is asymptotically equivalent to the process
%
%
\begin{equation}
\label{Udef2} \hat{U}_n^+ (t) \equiv n^{-1/2} \sum
_{i = 1}^{(X_n (0) - s_n (0))^+} 1_{\{\xi_i > t\}},\qquad t \ge0,
\end{equation}
where $\{\xi_i\}$ is a sequence of i.i.d. random variables, each
having a distribution with hazard rate $h_F (w (u))$
at time $u$. We only have asymptotic equivalence, because the
abandonment rate at time $u$ is actually $h_F (W_n (u))$ in system $n$.
However, we have $\|W_n - w\|_{\tau} \Rightarrow0$. Hence,\vspace*{1pt} for any
$\ep> 0$, we can bound the abandonment rate above by $h^{u, \ep}_F (w(u))$
and below by $h^{l, \ep}_F (w(u))$, where
%
%
\begin{eqnarray}
\label{hazBds} h^{u, \ep}_F \bigl(w(u)\bigr) & \equiv&\sup
_{-\ep\le s \le\ep}{\bigl\{h_F \bigl(w (u) + s\bigr)\bigr\}}\quad\mbox{and}
\nonumber\\[-8pt]\\[-8pt]
h^{l, \ep}_F \bigl(w(u)\bigr) & \equiv& \inf
_{-\ep\le s \le\ep}{\bigl\{ h_F \bigl(w (u) + s\bigr)\bigr\}}.\nonumber
\end{eqnarray}
By exploiting these bounds and the continuity of $f$, we see that we do
indeed asymptotically have the representation in (\ref{Udef2}).

Combining (\ref{InitialLim}), (\ref{Udef2}) and the Glivenko--Cantelli
theorem, we can conclude that
%
%
\begin{equation}
\label{abanLimitPlus} \hat{U}_n^+ (t) \Rightarrow\hat{X} (0)^+
F_w^c (t) \qquad\mbox{in } D\bigl([0,\tau]\bigr)
\end{equation}
as $n \ra\infty$. Essentially the same reasoning applies to
$A_{i,n,-} (t)$. Combining these two limits, we obtain
\[
\bigl\|\hat{X}_n - \bigl(\hat{X}^{*}_n +
\hat{X}_n (0) F^c_w (\cdot)\bigr)
\bigr\|_{\tau
} \Rightarrow0.
\]
Hence we have justified (\ref{asymEq2}).

\section{\texorpdfstring{Proofs of Theorem \protect\ref{thOther} and Corollary \protect\ref{corVar}}
{Proofs of Theorem 4.3 and Corollary 4.1}}\label{secCors}

\subsection{\texorpdfstring{Proof of Theorem \protect\ref{thOther}}
{Proof of Theorem 4.3}}

We have just proved that $\hat{X}_n \Rightarrow\hat{X}$ for the
general initial condition in (\ref{InitialLim}).
We now establish the remaining limits in (\ref{gXlims2}) for the other
related processes with initial condition (\ref{InitialLim}).

\subsubsection{\texorpdfstring{The processes $\hat{Q}_n$ and $\hat{B}_n$}
{The processes Q n and B n}}

Since\vspace*{1pt} $P(X_n (t) > s_n (t),   t_1 \le t \le t_2) \ra1$ as $n \ra
\infty$ for all
$t_1$ and $t_2$ with $0 < t_1< t_2 < \tau$,
we necessarily have $\|\hat{B}_n\|_{t_1, t_2} = \|\hat{X}_n - \hat
{Q}_n\|_{t_1, t_2} \Rightarrow0$
as $n \ra\infty$, so that $(\hat{X}_n, \hat{Q}_n, \hat{B}_n)
\Rightarrow(\hat{X}, \hat{X}, 0e)$
as claimed in $\D([t_1, t_2])$ for each $t_1$ and $t_2$ with $0 < t_1
< t_2 < \tau$,
which is equivalent to convergence in $\D((0, \tau))$.

However, the situation is different at the interval endpoints. In
particular, there is truncation at time $0$ for the
processes $\hat{Q}_n$ and $\hat{B}_n$.
We cannot extend the limit to the interval\vadjust{\goodbreak} $[0,\tau]$, or even $[0,
\tau)$, closed on the left, because the limit process
could have a discontinuity at $0$,
which would be ruled out in the definition of the space $\D$.
Indeed, because of the definition of the queue length as $Q_n (t)
\equiv(X_n (t) - s_n (0))^+$ and the number in service
as $B_n (t) \equiv X_n (t) \wedge s_n (0)$,
it is immediate that
%
%
\begin{equation}
\label{Qtrunc} \bigl(\hat{Q}_n (0),\hat{B}_n (0)\bigr)
\Rightarrow\bigl(\hat{X} (0)^+, \hat{X} (0)^-\bigr) \qquad\mbox{in }\RR^2
\mbox{ as } n
\ra\infty.
\end{equation}
These limits are of course not mean-zero random variables.

As a consequence, of (\ref{Qtrunc}),
if $P(\hat{X} (0) < 0) > 0$,
then there can be no FCLT for $\hat{Q}_n$ in $\D([0, \tau))$ because
$\hat{Q}$ would require a discontinuity at time $0$
to reflect the initial truncation of $X_n (0)$ to get $Q_n (0)$;
If $P(\hat{X} (0) > 0) > 0$,
then there is no FCLT for $\hat{B}_n$ in $\D([0, \tau))$
because $\hat{B}$ would require a discontinuity at time $0$
to reflect the initial truncation of $X_n (0)$ to get $B_n (0)$.

\subsubsection{\texorpdfstring{The abandonment process $A_n(t)$}
{The abandonment process A n(t)}}

We obtain the limits for $\bar{A}_n$ and $\hat{A}_n$ in (\ref{fLim})
and (\ref{gXlims2}) directly from
representation (\ref{aban}) and the established limits above. We see
that $\bar{A}_n \Rightarrow A$ in $\D$
as $n \ra\infty$ and
%
%
\begin{equation}
\label{Adiff} \sup_{0 \le t \le T}{\bigl\{\bigl|\hat{A}_n (t) -
\bigl(\hat{N}^{*}_n (t) - \hat {D}^{*}_n
(t) - \bigl(\hat{X}_n (t)\bigr) - \hat{X}_n (0)\bigr)\bigr|
\bigr\}} \Rightarrow0,
\end{equation}
so that $\hat{A}_n \Rightarrow\hat{A}$ in $\D$, jointly with the
other processes, for $\hat{A}$ in (\ref{Glims2}).

\subsubsection{The waiting times with the general initial conditions}

We will start by considering the PWT $V_n (t)$. We first consider time
$0$ for the FCLT-scaled process.
Note that the PWT $V_n (0)$ and the FCLT-scaled version $\hat{V}_n
(0)$ are $0$ if $X_n (0) \le0$, but not otherwise.
Hence, the general initial conditions in (\ref{InitialLim}) alters the
limit $\hat{V}$ at time $0$.
Since service times are exponential, service completion occurs
initially at rate
$(s_n (0) \wedge X_n (0)) \mu$, where $(s_n (0) \wedge X_n (0))/ n
\Rightarrow s(0)$. In addition, new service capacity initially
becomes available asymptotically at rate $n \dot{s} (0)$. Hence, the
scaled PWT at time $0$ is asymptotically equivalent to
%
%
\begin{equation}
\label{Vzero} \sqrt{n} \tilde{V}_n (0) \equiv n^{-1/2} \sum
_{i = 1}^{(X_n (0) -
s_n (0))^+} \zeta_i,
\end{equation}
where $\{\zeta_i\}$ is a sequence of i.i.d. exponential random
variables, each with rate $s(0) \mu+ \dot{s} (0) > 0$.
Hence, by the LLN,
%
%
\begin{equation}
\label{Vzero2} \hat{V}_n (0) = \sqrt{n} V_n (0)
\Rightarrow\frac{\hat{X}
(0)^+}{s(0)\mu+ \dot{s} (0)} \qquad\mbox{in }\RR\mbox{ as } n \ra\infty.
\end{equation}

We have a different situation for $t > 0$, because the number in system
becomes positive, of order $O(n)$ for $t>0$.
Since $P(X_n (t) > s_n (t),   t_1 \le t \le t_2) \Rightarrow1$ for
any $t_1$ and $t_2$ with $0 < t_1 < t_2 < \tau$,
now the service completion rate is asymptotically $s_n (t) \mu$ at
time $n$, for all $t$ in $[t_1, t_2]$ above.
As in (\ref{Udef}), we consider the remaining number from the initial
difference,
separating the positive and negative values.
Now, paralleling (\ref{Udef2}), we have $\hat{V}_n (t)$ asymptotically
equivalent to $\hat{V}_n^{*} (t) + \sqrt{n}\tilde{V}_n (t)$,
where $\tilde{V}_n (t) = \tilde{V}^{+}_n (t) - \tilde{V}^{-}_n (t)$ with
%
%
\begin{equation}
\label{Vpos} \sqrt{n}\tilde{V}^{+}_n (t) \equiv
n^{-1/2} \sum_{i = 1}^{U^{+}_n (t
+ V_n^{*} (t))}
\zeta_i,
\end{equation}
where $U^{+}_n (t)$ is defined in (\ref{Udef}) and $\{\zeta_i\}$ is a
sequence of i.i.d. exponential random variables, each with rate $s(t +
v(t))\mu+ \dot{s} (t + v(t))$,
and similarly for $\sqrt{n}\tilde{V}^{-}_n (t)$.
As a consequence, by the FWLLN,
%
%
\begin{equation}
\label{Vpos} \sqrt{n}\tilde{V}^{+}_n (t) \Rightarrow
\frac{\hat{X} (0)^+ F^c_w
(t + v (t))}{s (t + v(t))\mu+ \dot{s} (t+v(t))} \qquad\mbox{in }\D
((0, \tau)).
\end{equation}
Combining this result with the corresponding result for $\sqrt {n}\tilde{V}^{-}_n (t)$,
we have
%
%
\begin{eqnarray}
\label{Vgen2} \hat{V}_n (t) 
& \Rightarrow& \hat{V} (t) \equiv\hat{V}^{*} (t) +
\frac{\hat{X}
(0) F^c_w (t + v (t))}{s (t + v(t))\mu+ \dot{s} (t+v(t))} \qquad\mbox{in }\D
((0, \tau))
\nonumber
\\
& = & \frac{\hat{W}^{*} (t + v(t))}{1 - \dot{w} (t + v(t))} + \frac{\hat{X} (0) F^c_w (t + v (t))}{s (t + v(t))\mu+ \dot{s}
(t+v(t))}
\\
& = & \frac{\tilde{q} (t + v(t), w(t + v(t)))\hat{W}^{*} (t +
v(t)) + \hat{X} (0) F^c_w (t + v (t))}{s (t + v(t))\mu+ \dot{s}
(t+v(t))},\nonumber
\end{eqnarray}
where $\hat{V}^{*}$ has been determined already in Section~\ref{secFCLTother}, $F^c_w (t)$ is defined in (\ref{Glims})
and $b(t,0) = s(t) \mu+ \dot{s}(t) > 0$ by assumption in Section~\ref{secFluid}. The final formula in (\ref{Vgen2}) is equivalent to
the stated formula in (\ref{Glims2}) because $w(t + v(t)) = v(t)$ by
(\ref{vEq}).

We next use the equation $W_n (t) = V_n (t - W_n (t)) + O(1/n)$ to
develop a limit for $W_n (t)$.
Reasoning as in (\ref{v104}), we get $\hat{W}_n \Rightarrow\hat{W}$
in $D((0, \tau))$ with
%
%
\begin{equation}
\label{Vgen3} \hat{W} (t) = \bigl(1 - \dot{w}(t)\bigr)\hat{V} \bigl(t-v(t)\bigr)
\end{equation}
for $\hat{V}$ in (\ref{Vgen2}), from which the formulas given in
(\ref{Glims2}) follow directly.

\subsection{\texorpdfstring{Proof of Corollary \protect\ref{corVar}: The variance formulas}
{Proof of Corollary 4.1: The variance formulas}}\label{secVariances}

We obtain the complicated variance formulas for $\sigma^2_{\hat
{W}^{*}_i} (t)$ and $\sigma^2_{\hat{X}^{*}_i} (t)$
by applying the usual It\^{o} isometry for Brownian stochastic integrals,
using the representation in (\ref{Glims}).
The remaining variance formulas are elementary.

\section{\texorpdfstring{Proof of Theorem \protect\ref{thUL} for underloaded intervals}
{Proof of Theorem 5.1 for underloaded intervals}}\label{secULproofs}

In this section we prove Theorem \ref{thUL}.

\begin{pf*}{Proof of Theorem \ref{thUL}}
As indicated, this mostly is a direct application of the
infinite-server FWLLN\vadjust{\goodbreak} and FCLT in \cite{PW10}. This is true for $X_n$
because
%
%
\begin{equation}
\label{asymEq} \bigl\| \bar{X}_n - \bar{X}^{\infty}_n
\bigr\|_{\tau} \Rightarrow0 \quad\mbox{and}\quad\bigl\| \hat{X}_n -
\hat{X}^{\infty}_n\bigr\|_{\tau} \Rightarrow0,
\end{equation}
where $X_n^{\infty} (t)$
is the associated $G_t/M/\infty$ model with the identical arrival
process, the identical sequence of service times for successive
customers entering service
and the identical initial conditions, that is,
$X_n^{\infty} (0) \equiv X_n (0)$.
Thus we can apply many-server heavy-traffic (MSHT) limits established
for that model in \cite{PW10}; also see \cite{B67,KP97,RT09}.
(Previous references suffice here; the full force of \cite{PW10} is
only needed to treat the more general $G_t/\GI/\infty$ model associated
with OL intervals;
see Section~\ref{secOL}.)

However, to prove (\ref{asymEq}), we need to carefully consider what
happens in the neighborhood of each interval endpoint.
There is no trouble in between because there is no critical loading
except at the interval endpoints. That implies that the net flow out,
$D_n (t) - N_n (t) - n (s(t) - s(0))$, is positive of order $O(n)$
over any interval $[t_1,t_2]$ for $0 < t_1 < t_2 < \tau$, no matter
how small. Thus, $P(\sup_{t_1 \le s \le t_2}{\{X_n (s) - s_n (s)\}} <
0) \ra1$
as $n \ra\infty$ for $0 < t_1 < t_2 < \tau$.\vspace*{1pt}

However, it is possible that $X_n (0) > s_n (0)$ and/or $X_n (\tau) >
s_n (\tau)$. Consider the left endpoint. If $X_n (0) > s_n (0)$, then
the systems
$X_n$ and $X_n^{\infty}$ are not stochastically identical over $[0,t]$
for $t > 0$.
We do have $X_n^{\infty} (0) = X_n (0)$ by definition, but if $X_n (0)
> s_n (0)$, then $X_n (0) - s_n (0)$ customers are waiting in queue
instead of being served.
However, asymptotically, the difference at time $0$ is $\sqrt{n} \hat
{X}(0)^{+} = O(\sqrt{n})$.
Only this portion of the initial number of customers will receive
different treatment.

Since $\delta(t) \equiv s(t) - X(t)$ is differentiable with derivative
$\dot{\delta} (0) > 0$, the initial difference of order $O(\sqrt{n})$
will dissipated over a time interval of
order $O(1/\sqrt{n})$. The constant departure rates (by service versus
abandonment) of these $O(\sqrt{n})$
customers will differ during that short time interval. Thus, $\|X_n -
X_n^{\infty}\|_t$ is of order $O(\sqrt{n}) \times O(1/\sqrt{n}) =
O(1)$ as $n \ra\infty$.
Hence, this difference is asymptotically negligible after scaling.
To support this conclusion, note that the hazard rate of the
abandonment is bounded above, implying that only a negligible number of
customers in the queue will
abandon in the initial interval of length $O(1/\sqrt{n})$.

Essentially the same argument applies at the right endpoint $\tau$.
Thus, we do indeed have $\| \bar{X}_n - \bar{X}^{\infty}_n\|_{\tau}
\Rightarrow0$ and $\| \hat{X}_n - \hat{X}^{\infty}_n\|_{\tau}
\Rightarrow0$,
as claimed in (\ref{asymEq}). Then the results for the $G_t/M/\infty$
model follow from \cite{PW10}. A key step there is to treat the new
arrivals differently from the
customers initially in the system. The customers initially in the
system are treated in Section~5 of \cite{PW10}; they lead to the limit
processes
$X_z$ and $\hat{X}_z$.

However, truncation at the endpoints $0$ and $\tau$ do alter the
processes $B_n$ and $Q_n$ more significantly.
Since we can have $\bar{X}_n (0) \neq s(0)$,
and/or
$\bar{X}_n (\tau) \neq s(\tau)$ for all~$n$, there can be
truncation at the times $0$ and $\tau$. Thus we can have
$\bar{B}_n (0) = s(0) \neq \bar{X}_n (0)$ and/or $\bar{B}_n (\tau
) = s(\tau) \neq \bar{X}_n (\tau)$. However, there is no problem
for the fluid limits.
Since $\bar{X} (0) \le s(0)$ and $\hat{X}_n (0) \Rightarrow\hat{X}
(0)$, necessarily
$\|\bar{X}_n (0) - \bar{B}_n (0) \| = O(1/\sqrt{n}) = o(1)$, so that
(\ref{BfluidLim}) follows from Theorem 11.4.7 of \cite{W02}.
The same reasoning can be applied at the right endpoint $\tau$.

In contrast, the truncation affects the FCLTs for $\hat{Q}_n$ and
$\hat{B}_n$ when $X(0) = s_n (0)$
Since $P(X_n (t) < s_n (t),  0< t_1 \le t \le t_2 < \tau) \ra1$
as $n \ra\infty$,
we necessarily have $\|\hat{Q}_n\|_{t_1, t_2} = \|\hat{X}_n - \hat
{B}_n\|_{t_1, t_2} \Rightarrow0$
as $n \ra\infty$, so that $(\hat{X}_n, \hat{B}_n, \hat{Q}_n)
\Rightarrow(\hat{X}, \hat{X}, 0e)$
as claimed. We cannot extend the limit to the closed interval $[0, \tau
]$ because the limit process could have a discontinuity at $0$,
which would be ruled out.
If $P(\hat{X} (0) < 0) > 0$,
then there can be no FCLT for $\hat{Q}_n$ in $\D([0, \tau))$ because
$\hat{Q}$ would require a discontinuity at time $0$
to reflect the initial truncation of $X_n (0)$ to get $Q_n (0)$;
If $P(\hat{X} (0) > 0) > 0$,
then there is no FCLT for $\hat{B}_n$ in $\D([0, \tau))$
because $\hat{B}$ would require a discontinuity at time $0$
to reflect the initial truncation of $X_n (0)$ to get $B_n (0)$.
There also could be further truncation at the right endpoint $\tau$,
so we only state
the limit for $(\hat{B}_n, \hat{Q}_n)$ in $\D([0, \tau))$.
\end{pf*}

Extending Theorem \ref{thUL} to the more general $G_t/\GI/s_t+\GI$ model
is more difficult, because the limit
for $\bar{X}_{z,n}$ involving the initial customers would be more
complicated because it would depend on the
ages of all the service times in process. We have exploited the
exponential assumption to avoid that difficulty.

\section{\texorpdfstring{Comparison with simulation: An $M_{t}/M/s_{t}+H_{2}$ example}
{Comparison with simulation: An M t/M/s t + H 2 example}}\label{secCompare}

To provide practical confirmation of the theorems proved in earlier
sections, we now report the results of a simulation experiment.
We consider an $M_t/M/s+H_2$ queueing model with a sinusoidal arrival
rate function that makes the system alternate between OL and UL intervals.
Specifically, the model parameters are: arrival rate function
$\lambda_{n}(t) = n\lambda(t)$, $\lambda(t) = 1 + 0.6 \sin(t)$,
mean service time $1/\mu= 1$, mean patience
$1/\theta= 2$ and a fixed number of servers $s_{n}(t) = n s$, $s =
1$. We let the service distribution be exponential and the patience
distribution be
a two-phase hyperexponential ($H_2$) with probability density function (p.d.f.)
\[
f(x) = p\cdot\theta_{1}e^{-\theta_{1}x} + (1-p)\cdot\theta
_{2}e^{-\theta_{2}x},\qquad x\geq0,
\]
with parameters $p=0.5(1-\sqrt{0.6})$, $\theta_{1}=2p\theta$ and
$\theta_{2}=2(1-p)\theta$, which produces squared coefficient of
variation (variance divided by the square of the mean)
\mbox{$c^{2} = 4$}.

\begin{figure}

\includegraphics{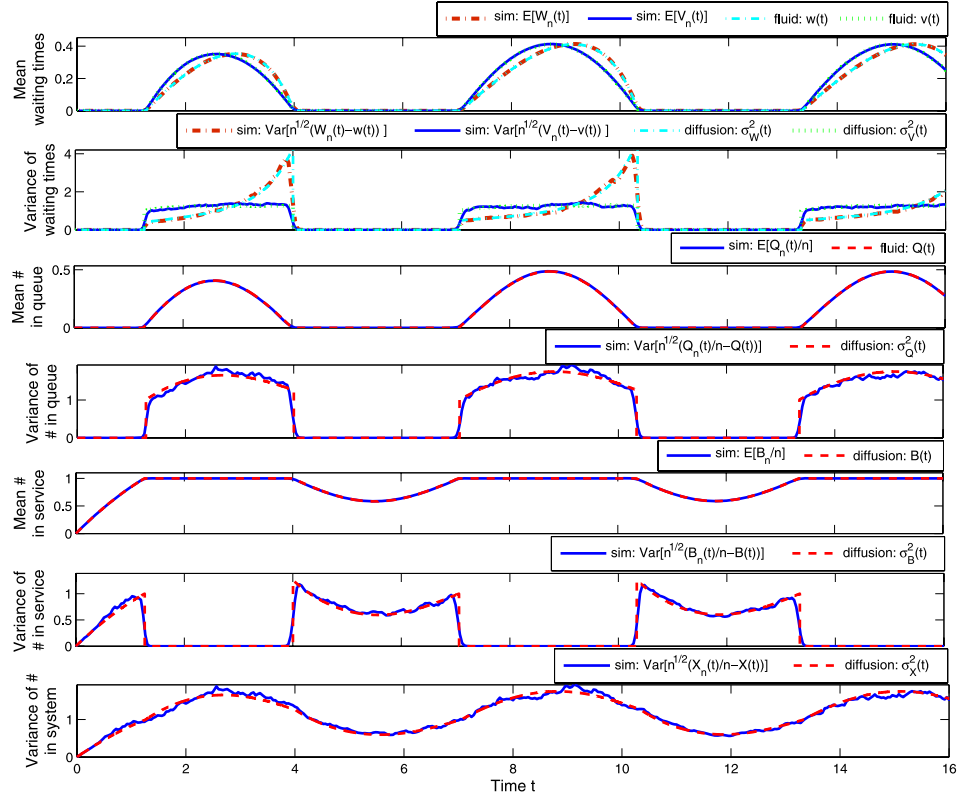}

\caption{Comparison of the limiting means (fluid limits) and variances
of the Gaussian limits to simulation estimates
of the corresponding scaled queueing processes for the $M_t/M/s+H_2$
model starting empty for the case $n = 2000$ based on $500$
independent replications: \textup{(i)} the boundary and potential waiting
times, $w(t)$ and $v(t)$,
\textup{(ii)} the variances of the two waiting times,
\textup{(iii)} mean number in queue, $Q(t)$, \textup{(iv)} the variance of the
number in queue, \textup{(v)} mean number in service, $B(t)$,
\textup{(vi)}~variance of number in service and \textup{(vii)} variance of the
total number in the system, $X(t)$.}
\label{H2abfig}
\end{figure}

To verify accuracy of the formulas, we estimate the mean and variance
of the scaled queueing processes for very large $n$, in particular, for
$n = 2000$.
We obtain these estimates from $500$ independent replications of a
simulation of the queueing system.
Figure~\ref{H2abfig} shows plots of several key performance
functions for the limiting fluid and diffusion processes for $0\leq
t\leq T\equiv16$,
starting out empty (see dashed lines):
(i) \textit{fluid} head-of-line and the potential waiting times $w(t)$ and $v(t)$,
(ii) variance of the \textit{diffusion} waiting times $\sigma^{2}_{\hat
{W}}(t)$ and $\sigma^{2}_{\hat{V}}(t)$,
(iii) \textit{fluid} number of customers in queue, in service $Q(t)$ and $B(t)$,
(iv) variance of the \textit{diffusion} number of customers in queue, in
service, and in the system $\sigma^{2}_{\hat{Q}}(t)$, $\sigma
^{2}_{\hat{B}}(t)$, and $\sigma^{2}_{\hat{X}}(t)$.

We compare these performance functions for the limit processes to
estimates of them for the corresponding scaled queueing processes.
In Figure~\ref{H2abfig} we also plot the corresponding performance
functions under the LLN and CLT scaling (see solid lines):
(i) mean of the LLN-scaled head-of-line and the potential waiting times
$E[\bar{W}_{n}(t)]$ and $E[\bar{V}_{n}(t)]$,
(ii) variance of the CLT-scaled waiting times $\operatorname{Var}(\hat{W}_{n}(t))$
and $\operatorname{Var}(\hat{V}_{n}(t))$,
(iii) mean of the LLN-scaled number of customers in queue and in
service $E[\bar{Q}_{n}(t)]$ and $E[\bar{B}_{n}(t)]$,
(iv) variance of the CLT-scaled number of customers in queue, in
service, and in the system $\operatorname{Var}(\hat{Q}_{n}(t))$, $\operatorname{Var}(\hat
{B}_{n}(t))$ and $\operatorname{Var}(\hat{X}_{n}(t))$.
Figure~\ref{H2abfig} shows that the simulation estimates for the
$M_t/M/s+H_2$ queueing model agree closely with the fluid and diffusion
performance.

\begin{figure}

\includegraphics{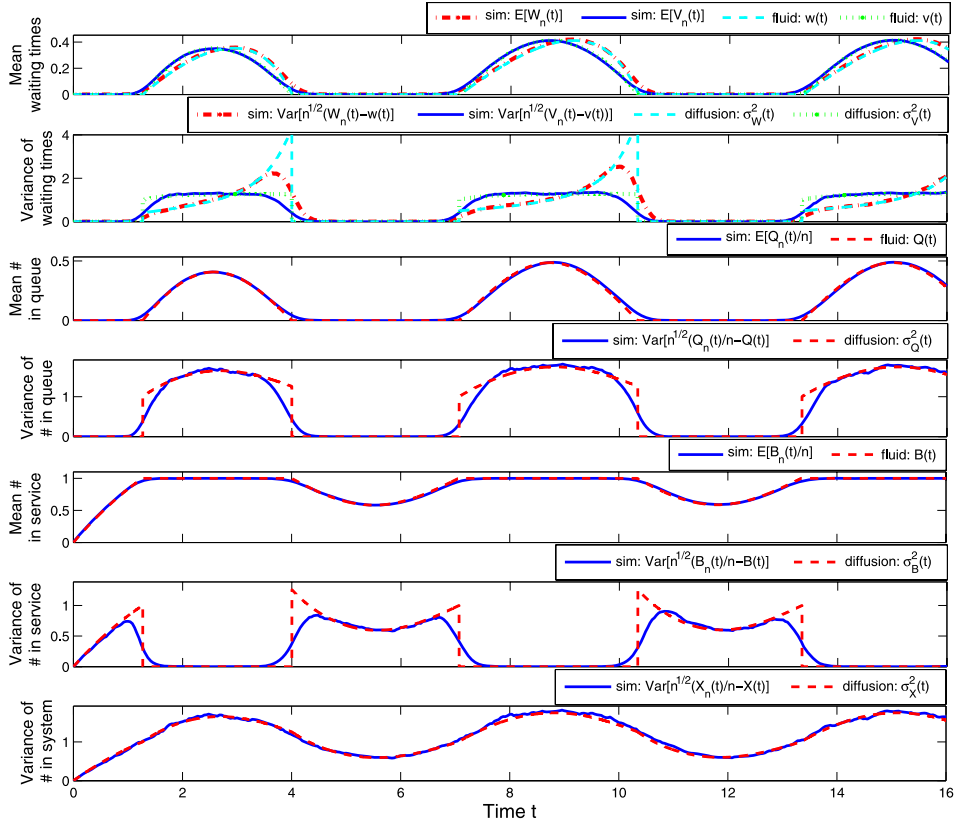}

\caption{Comparison of the limiting means (fluid limits) and variances
of the Gaussian limits to simulation estimates
of the corresponding scaled queueing processes for the $M_t/M/s+H_2$
model starting empty for the case $n = 100$ based on $2000$
independent replications: \textup{(i)} the boundary and potential waiting
times, $w(t)$ and $v(t)$,
\textup{(ii)} the variances of the two waiting times,
\textup{(iii)} mean number in queue, $Q(t)$, \textup{(iv)} the variance of the
number in queue, \textup{(v)} mean number in service, $B(t)$,
\textup{(vi)}~variance of number in service and \textup{(vii)} variance of the
total number in the system, $X(t)$.}
\label{H2abfig100}
\end{figure}

This experiment provides an engineering verification for the limit
theorems (as $n\rightarrow\infty$).
The approximation is not nearly as good when $n$ is small, for example,
when $n = 100$, as shown in Figure~\ref{H2abfig100}. The
approximation still performs well in the interior of UL and OL
intervals but relatively poorly in the neighborhood of the switching
points (the real variances are continuous functions while the
approximating formulas are jump functions). Furthermore, we find the
approximation becomes even worse for smaller systems, for example, when
$n=20$. Thus, we develop and study refined engineering approximations,
drawing on (\ref{QBapprox}), in \cite{LW11b}.\

\section{\texorpdfstring{Refined scaling with additional $O(\sqrt{n})$ terms}
{Refined scaling with additional O(square of n) terms}}
\label{secRefined}

For refined approximations and controls, we may want to generalize the
sequence of $G_t/M/s_t+\GI$ queueing models specified in Section~\ref{secSeq} by considering
arrival rates
$\lambda_n (t) \equiv n \lambda(t) + \sqrt{n} \lambda_g (t)$
and staffing functions
$s_n (t) \equiv\lceil ns(t) + \sqrt{n} s_g (t)\rceil$,
having extra $\sqrt{n}$ terms, where $\lambda_g (t)$ and $s_g (t)$
are additional smooth deterministic functions (with subscript $g$ for
Gaussian scale).
We now briefly indicate how the results above extend to this case.

First, the limit processes
in the FCLT for the arrival process and the departure process in (\ref
{arriveFCLT}) and (\ref{twoFCLT})
should have the respective extra terms
%
%
\begin{equation}
\label{refineAD} \Lambda_g (t) \equiv\int_{0}^{t}
\lambda_g (s) \,ds \quad\mbox{and}\quad D_g (t) \equiv\mu\int
_{0}^{t} s_g (s) \,ds.
\end{equation}
These changes lead to deterministic modifications of other expressions.

For each OL interval,
we add the term $Z_{1,g} (t,y) \equiv\int_{t - y}^{t} F^c (t-s)
\lambda_g (s) \,ds$ to $\hat{Z}_{1} (t,y)$
in (\ref{tildeLimits}); we add the term $Q_{1,g} (t) \equiv\int_{t -
w(t)}^{t} F^c (t-s) \lambda_g (s) \,ds$
to $\hat{Q}_{1} (t)$ in Section~\ref{secQlimPf}; we add the tem $D_g
(t)$ to $\hat{D} (t)$ in (\ref{Elimit});
and we add the term $I_{1,g} (t, \ep) \equiv\int_{t-w(t)}^{t+\varepsilon-w(t+\varepsilon)}F^c(t-s) \lambda_g (s) \,ds$
to $\hat{I}_1 (t, \ep)$ in (\ref{hatIn1}).

Those changes lead
to changes in the critical SDE for the limit process $\hat{W} (t)$
developed in Section~\ref{secWcharPf}. Extra terms
$D_g (t+ \ep) - D_g (t)$ appear on the left and
$\int_{t-w(t)}^{t+\varepsilon-w(t+\varepsilon)}F^c(t-s) \lambda_g (s) \,ds$ on the right in (\ref{preSDE}),
which in turn contribute a term $-z(t) \,dt$ to the right-hand side of
the SDE in (\ref{sdePf}) and (\ref{sde}), where
%
%
\begin{equation}
\label{zDef} z(t) \equiv\frac{s_g (t) \mu+ \lambda_g (t - w(t))+ \dot{s}_g (t)}{q(t,w(t))}.
\end{equation}
This leads to an extra deterministic term $W_g (t) \equiv- \int_{0}^{t}H(t,u) z(u)\,d u$ on the right-hand side of the expression for
$\hat{W} (t)$ given in (\ref{SDEsol}), which is $\hat{W}^{*} (t)$ in
Theorem~\ref{thFCLT}.

From (\ref{Q1}) and (\ref{Qlim}), we see that
those changes above lead to the addition of $Q_{1,g} (t)$ above to
$\hat{X}^{*}_1 (t)$ and the addition of $q(t, w(t)) W_g (t)$ to $\hat
{X}^{*}_3 (t)$
in Theorem~\ref{thFCLT}.

There are corresponding changes for each UL interval. Due to the
revised arrival and departure FCLTs,
the term $u(t) \,dt$ is added to the right-hand side of the SDE in
(\ref{sdeUL}), where $u(t) \equiv\lambda_g (t)$.

The changes above lead to modifications of the limits in the FCLTs,
but not the FWLLNs. The limit processes are still Gaussian processes.
These deterministic changes alter the mean values of the Gaussian
limits, but do not affect the variances.

\section*{Acknowledgment}

We thank Jiheng Zhang for constructive comments.



\printaddresses

\end{document}